\documentclass[12 pt, onecolumn]{IEEEtran}
\usepackage{amsmath,amssymb,euscript,yfonts,psfrag,latexsym,graphicx}
\usepackage{bbm,color,amstext,wasysym,subfig,flushend,parskip,balance}
\graphicspath{{./},{./figures/}}

\newtheorem{thm}{Theorem}

\newtheorem{remark}[thm]{Remark}

\newcommand{\mR}{{\mathbb R}}

\newcommand{\mD}{{\mathbb D}}

\newcommand{\cD}{{\mathcal D}}
\newcommand{\cE}{{\mathcal E}}
\newcommand{\cF}{{\mathcal F}}

\newcommand{\cH}{{\mathcal H}}

\newcommand{\cS}{{\mathcal S}}

\newcommand{\cV}{{\mathcal V}}

\newcommand{\cL}{{\mathcal L}}
\newcommand{\cW}{{\mathcal W}}
\newcommand{\bw}{{\mathbf w}}
\newcommand{\diag}{\operatorname{diag}}

\newcommand{\trace}{\operatorname{tr}}

\definecolor{grey}{rgb}{0.6,0.6,0.6}
\definecolor{lightgray}{rgb}{0.97,.99,0.99}

\begin{document}
\title{An Efficient Algorithm for Matrix-Valued and Vector-Valued Optimal Mass Transport}
\author{Yongxin Chen, Eldad Haber, Kaoru Yamamoto, Tryphon T. Georgiou, and Allen Tannenbaum
\thanks{Y.\ Chen is with the Department of Medical Physics, Memorial Sloan Kettering Cancer Center, NY; email: chen2468@umn.edu}
\thanks{E.\ Haber is with the Department of Mathematics, University British Columbia, Vancouver, Canada; email: haber@math.ubc.ca}
\thanks{K.\ Yamamoto is with the Department of Electrical Engineering, Lund University, Sweden; email: k.yamamoto@ieee.org}
\thanks{T.\ T. Georgiou is with the Department of Mechanical and Aerospace Engineering, University of California, Irvine, CA; email: tryphon@uci.edu}
\thanks{A.\ Tannenbaum is with the Departments of Computer Science and Applied Mathematics \& Statistics, Stony Brook University, NY; email: allen.tannenbaum@stonybrook.edu}}

\maketitle

\begin{abstract}
We present an efficient algorithm for recent generalizations of optimal mass transport theory to matrix-valued and vector-valued densities.
These generalizations lead to several applications including diffusion tensor imaging, color images processing, and multi-modality imaging.
The algorithm is based on sequential quadratic programming (SQP). By approximating the Hessian of the cost and solving each iteration
in an inexact manner, we are able to solve each iteration with relatively low cost while still maintaining a fast convergent rate.
The core of the algorithm is solving a weighted Poisson equation, where different efficient preconditioners may be employed.
We utilize incomplete Cholesky factorization, which yields an efficient and straightforward solver for our problem.
Several illustrative examples are presented for both the matrix and vector-valued
cases.

\end{abstract}

\section{Introduction}
The theory of optimal mass transport (OMT) \cite{Rac98,Vil03,AmbGigSav06} has proven its power and usefulness in both theory and applications.
The theory part has been developed through a sequence of elegant papers, and the research is still going strong;
see \cite{Mon81,Kan42,Bre91,GanMcc96,Mcc97,JorKinOtt98,BenBre00,OttVil00} and the references therein.
On the other hand, during the past decade, the need for applications has engendered the fast development of efficient
algorithms for OMT \cite{AngHakTan03,Cut13,BenFroObe14,HabHor15,BenCarCut15,CheGeoPav15a,LiYinOsh16,LiRyuOsh17}.
Recently, the OMT theory has been extended to study matrix \cite{CarMaa16,CheGeoTan16,MitMie16} and vector-valued densities \cite{CheGeoTan17}.

The mathematical approach to matrix optimal mass transport in \cite{CarMaa16,CheGeoTan16,MitMie16} is based on the seminal work of Benamou-Brenier \cite{BenBre00}, where
optimal mass transport with quadratic cost is recast as the problem of minimizing kinetic energy (i.e., an action integral) subject to a continuity equation. In the matrix case, one needs to develop a non-commutative counterpart to scalar optimal transport where probability distributions are replaced by density matrices $\rho$ (Hermitian positive-definite with unit trace) and where ``transport'' corresponds to a flow on the space of such matrices that minimizes a corresponding action integral. The work is motivated by a plethora of applications including spectral analysis of vector-valued time-series, which may encode different modalities (e.g., frequency, color, polarization) across a distributed array of sensors \cite{NinGeoTan15}. The associated power spectra are matrix-valued and hence there is a need for suitable metrics that quantify distances and provide tools to
average and interpolate spectra.  The generalization of the Benamou-Brenier theory is founded upon concepts from quantum mechanics, and allows us to
formulate a continuity equation for matrix-flows, and then derive a Wasserstein distance between density matrices and matrix-valued distributions.
Similar remarks apply to the vector-valued case in which one must also invoke some ideas from graph theory in formulating our generalization of scalar-valued densities.
See \cite{CheGeoTan17} for all the details.

In this paper, we focus on algorithms for the numerical solution of the optimal matrix-valued mass transport problems introduced in \cite{CarMaa16,CheGeoTan16,MitMie16},
and the vector-valued case formulated in \cite{CheGeoTan17}. In \cite{CheGeoTan16,CheGeoTan17}, both problems are reformulated as convex optimization problems.
We adopt an inexact sequential quadratic programming (SQP) method \cite{SteHab17,ByrCurNoc08,NocWri06} to tackle such convex optimization problems.
Similar methods have been applied to scalar optimal mass transport \cite{HabHor15}.

The remainder of  this paper is summarized as follows. Section~\ref{sec:background} is a brief introduction to the matrix-valued optimal transport theory.
We develop the corresponding algorithm in Section~\ref{sec:algorithm}, and then the algorithm for vector-valued optimal transport is described in Section~\ref{sec:vector}.
We conclude with several examples to demonstrate our algorithm in Section~\ref{sec:numerical}.

\section{Matrix-valued optimal mass transport} \label{sec:background}

In this section, we sketch the approach \cite{CheGeoTan16} for which the convex optimization algorithm given in the present note was formulated.
As noted above, similar approaches to matrix-valued OMT were formulated independently in \cite{CarMaa16,MitMie16}.

\subsection{Gradient on space of Hermitian matrices}\label{sec:quantumgrad}
Denote by $\cH$ and $\cS$ the set of  $n\times n$ Hermitian and skew-Hermitian matrices, respectively. We will assume that all of our matrices
are of fixed size $n\times n$. Next, we denote the space of block-column vectors consisting of $N$ elements in $\cS$ and $\cH$ as $\cS^N$ and $\cH^N$, respectively.
We also let $\cH_+$ and $\cH_{++}$ denote the cones of nonnegative and positive-definite matrices, respectively,
and we use the standard notion of inner product, namely
	\[
		\langle X,Y\rangle=\trace(X^*Y),
	\]
for both $\cH$ and $\cS$.
For $X, Y\in \cH^N$ ($\cS^N$),
	\[
		\langle X, Y\rangle=\sum_{k=1}^N \trace(X_k^*Y_k).
	\]
Given $X=[X_1^*,\cdots,X_N^*]^* \in \cH^N$ ($\cS^N$), $Y\in \cH$ ($\cS$), set
	\[
		XY=\left[\begin{array}{c}
		X_1\\
		\vdots \\
		X_N
		\end{array}\right]Y
		:=
		\left[\begin{array}{c}
		X_1Y\\
		\vdots \\
		X_N Y
		\end{array}\right],
	\]
	\[
		YX=Y\left[\begin{array}{c}
		X_1\\
		\vdots \\
		X_N
		\end{array}\right]
		:=
		\left[\begin{array}{c}
		YX_1\\
		\vdots \\
		YX_N
		\end{array}\right],
	\]
and	
	\[
		\bar X = \left[\begin{array}{c}
		X_1^*\\
		\vdots \\
		X_N^*
		\end{array}\right].
	\]

For a given $L\in\cH^N$ we define
	\begin{equation}\label{eq:gradient}
		\nabla_L: \cH \rightarrow {\cS}^N, ~~X \mapsto
		\left[ \begin{array}{c}
		L_1 X-XL_1\\
		\vdots \\
		L_N X-X L_N
		\end{array}\right]
	\end{equation}
to be the \textbf{\emph{gradient operator}}.
By analogy with the ordinary multivariable calculus, we refer to its dual with respect to the Hilbert-Schmidt inner product as the (negative)
\textbf{\emph{divergence operator}}, and this is
	\begin{equation}\label{eq:divergence}
		\nabla_L^*: {\cS}^N \rightarrow \cH,~~Y=
		\left[ \begin{array}{c}
		Y_1\\
		\vdots \\
		Y_N
		\end{array}\right]
		\mapsto
		\sum_k^N L_k Y_k-Y_k L_k,
	\end{equation}
i.e., $\nabla_L^*$ is defined by means of the identity
	\[
		\langle \nabla_L X , Y\rangle =\langle X , \nabla_L^* Y\rangle.
	\]
A standing assumption throughout, is that the null space of $\nabla_L$, denoted by ${\rm ker}(\nabla_L)$, contains only scalar multiples of the
identity matrix. In this note, we use one such basis generated by the following $N=2$ components:
	\[
		L_1=\left[\begin{matrix}1 & 1 &\cdots &1\\1& 0 &\cdots & 0\\\vdots &\vdots
		&\ddots & \vdots\\1& 0 & \cdots& 0\end{matrix}\right],
		\quad
		L_2=\diag ([1,\,2,\,\ldots,\,n-1,\,0]).
	\]
\subsection{Matrix-valued Optimal mass transport}
We next sketch the formulation for matrix-valued optimal mass transport proposed in \cite{CheGeoTan16}. Given a convex compact set $E\in \mR^m$, denote
	\[
		\cD = \{\rho(\cdot) \in \cH_+ ~\mid~ \int_{E} \trace(\rho(x)) dx =1\},
	\]
and $\cD_+$ the interior of $\cD$. Let $\rho^0,\rho^1\in \cD_+$ be two matrix-valued densities defined on $E$ with positive values.
A dynamic formulation of matrix-valued optimal mass transport between
these two given marginals is \cite{CheGeoTan16},
	 \begin{subequations}\label{eq:matrixomt}
 	\begin{eqnarray}\label{eq:matrixomt1}
	&&\min_{\rho\in \cD_+, w\in {\cH}^m, v\in \cS^N} \int_0^1\int_{E}
	\left \{\trace(\rho w^*w)+{\gamma}\trace(\rho v^*v)\right\}dxdt,\\
	&&\frac{\partial \rho}{\partial t}+\frac{1}{2}\nabla_x\cdot(w\rho+\rho w)-\frac{1}{2} \nabla_L^* (v\rho+\rho v)=0, \label{eq:matrixomt2}\\
	&& \rho(0,\cdot)=\rho^0, ~~\rho(1,\cdot)=\rho^1  \label{eq:matrixomt3}
	\end{eqnarray}
	\end{subequations}
with $\nabla_x \cdot$ being the standard divergence operator in $\mR^m$. By defining $p=w\rho, u=v\rho$, the above can be cast as a convex optimization problem
	 \begin{subequations}\label{eq:matrixomtcvx}
 	\begin{eqnarray}\label{eq:matrixomtcvx1}
	&&\min_{\rho, p, u} \int_0^1\int_{E}
	\left \{\trace(p\rho^{-1}p^*)+{\gamma}\trace(u\rho^{-1}u^*)\right\}dxdt\\
	&&\frac{\partial \rho}{\partial t}+\frac{1}{2}\nabla_x\cdot(p+\bar p)-\frac{1}{2} \nabla_L^*
	(u-\bar u)=0, \label{eq:matrixomtcvx2}\\
	&& \rho(0,\cdot)=\rho^0, ~~\rho(1,\cdot)=\rho^1  \label{eq:matrixomtcvx3}.
	\end{eqnarray}
	\end{subequations}
We remark that $(p+\bar p)/2\in \cH^m$ and $(u-\bar u)/2\in\cS^N$, which is consist with the domain of $\nabla_L^*$.
For the sake of brevity, the set $E$ is taken to be the unit cube $[0,1]^m$.

\section{Discretization and algorithm: matrix-valued case} \label{sec:algorithm}

We follow closely the algorithm developed in \cite{HabHor15} for scalar optimal mass transport problems.
We restrict ourselves to the real-valued case, that is, $\cH$ and $\cS$ denote symmetric and skew-symmetric matrices, respectively.
In order to highlight the key parts of our methodology, we first consider the discretization in 1D case, i.e., $m=1$.
In particular, we take $E=[0,1]$. The algorithm extends almost verbatim to the higher dimensional setting as we will see in Section~\ref{sec:2D}.

We discretize the space-time domain $[0,\,1]\times[0,\,1]$ into $n_x\times n_t$ rectangular cells. Denote $\Omega_{ij}, 1\le i\le n_x, 1\le j\le n_t$ as the $(i,j)$ box.
We use a staggered grid to discretize $p$ and $\rho$. {\bf \em The variable $u$ is, however, valued at the centers of the cells $\{\Omega_{ij}\}$}. More specifically,
	\[
		p=(p_{i+\frac{1}{2},j}), ~0\le i\le n_x, ~1\le j\le n_t
	\]
	\[
		\rho=(\rho_{i,j+\frac{1}{2}}),~1\le i\le n_x,~0\le j\le n_t
	\]
	\[
		u=(u_{i,j}),~1\le i \le n_x, ~1\le j\le n_t.
	\]
Note the boundary values are
	\[
		p_{\frac{1}{2},j}=0, ~p_{n_x+\frac{1}{2},j}=0, ~1\le j\le n_t
	\]
and
	\[
		\rho_{i,\frac{1}{2}}=\rho^0_i, ~\rho_{i,n_t+\frac{1}{2}}=\rho^1_i, ~1\le i\le n_x.
	\]
We exclude the boundary values from the variables and denote
	\[
		p=(p_{i+\frac{1}{2},j}), ~1\le i\le n_x-1, ~1\le j\le n_t
	\]
	\[
		\rho=(\rho_{i,j+\frac{1}{2}}),~1\le i\le n_x,~1\le j\le n_t-1.
	\]
	
\subsection{Continuity equation}

We use the above discretizing scheme, together with the boundary conditions to rewrite the continuity equation \eqref{eq:matrixomtcvx2} as
	\begin{equation}\label{eq:discretecontinuity}
	D_1 p+D_2 \rho+D_3 u=b.
	\end{equation}
Here the linear operators $D_1, D_2, D_3$ are defined as
	\[
		(D_1 p)_{i,j}=
		\begin{cases}
		\frac{1}{2}(p_{i+\frac{1}{2},j}+p_{i+\frac{1}{2},j}^*
		-p_{i-\frac{1}{2},j}-p_{i-\frac{1}{2},j}^*)/h_x,  & \mbox{if} ~~2\le i \le n_x-1,
		\\
		\frac{1}{2} (p_{\frac{3}{2},j}+p_{\frac{3}{2},j}^*)/h_x, & \mbox{if} ~~i=1,
		\\
		-\frac{1}{2}(p_{n_x-\frac{1}{2},j}+p_{n_x-\frac{1}{2},j}^*)/h_x, & \mbox{if} ~~i=n_x,
		\end{cases}
	\]
	\[
		(D_2 \rho)_{i,j}=
		\begin{cases}
		(\rho_{i,j+\frac{1}{2}}-\rho_{i,j-\frac{1}{2}})/h_t,  & \mbox{if} ~~2\le j \le n_t-1,
		\\
		 \rho_{i,\frac{3}{2}}/h_t, & \mbox{if} ~~j=1,
		\\
		- \rho_{i,n_t-\frac{1}{2}}/h_t, & \mbox{if} ~~j=n_t,
		\end{cases}
	\]
	\[
		(D_3 u)_{i,j}=-\frac{1}{2}\nabla_L^*(u_{i,j}-\bar u_{i,j}), ~1\le i \le n_x, ~1\le j\le n_t.
	\]
The parameter $b$ carries the information of the boundary values $\rho^0$ and $\rho^1$. More specifically,
	\[
		b_{i,j}=
		\begin{cases}
		\rho^0_i/h_t & \mbox{if}~~j=1,\\
		-\rho^1_i/h_t & \mbox{if}~~j=n_t,\\
		0 & \mbox{otherwise}.
		\end{cases}
	\]

\subsection{Discretizing the cost function}
We use a combination of a midpoint and a trapezoidal methods to discretize the cost function. On the volume $\Omega_{ij}$ we have
	\begin{eqnarray*}
		\int_{\Omega_{ij}} \left \{\trace(p\rho^{-1}p^*)+{\gamma}\trace(u\rho^{-1}u^*)\right\}
		&\approx& \frac{h_xh_t}{4}\trace((p_{i-\frac{1}{2},j}^*p_{i-\frac{1}{2},j}+
		p_{i+\frac{1}{2},j}^*p_{i+\frac{1}{2},j})(\rho_{i,j-\frac{1}{2}}^{-1}+\rho_{i,j+\frac{1}{2}}^{-1}))
		\\&&+\frac{\gamma h_xh_t}{2}\trace(u_{i,j}^*u_{i,j}(\rho_{i,j-\frac{1}{2}}^{-1}+
		\rho_{i,j+\frac{1}{2}}^{-1})).
	\end{eqnarray*}
Let $A_1$ be the averaging operator over the spatial domain and $A_2$ be the averaging operator over the time domain (one needs to be careful about the boundaries).
Then the cost function  \eqref{eq:matrixomtcvx1} may be approximated by
	\begin{equation}\label{eq:costdiscrete}
		\left\langle A_1 (p^*\circ p), A_2 (\rho^{-1})+a\right\rangle h_xh_t+
		\left\langle u^*\circ u, A_2 (\rho^{-1})+a\right\rangle\gamma h_xh_t,
	\end{equation}
where $a\ge 0$ depends only on the boundary values $\rho^0$ and $\rho^1$. The inverse operator and the multiplication operator $\circ$ are applied block-wise.
The expressions for $A_1, A_2, a$ are
	\[
		(A_1 (p^*\circ p))_{i,j}=
		\begin{cases}
		\frac{1}{2}(p^*_{i-\frac{1}{2},j}p_{i-\frac{1}{2},j}
		+p^*_{i+\frac{1}{2},j}p_{i+\frac{1}{2},j}),  & \mbox{if} ~~2\le i \le n_x-1,
		\\
		\frac{1}{2}p^*_{\frac{3}{2},j}p_{\frac{3}{2},j}, & \mbox{if} ~~i=1,
		\\
		\frac{1}{2}p^*_{n_x-\frac{1}{2},j}p_{n_x-\frac{1}{2},j}, & \mbox{if} ~~i=n_x,
		\end{cases}
	\]
	\[
		(A_2 (\rho^{-1}))_{i,j}=
		\begin{cases}
		\frac{1}{2}(\rho^{-1}_{i,j-\frac{1}{2}}
		+\rho^{-1}_{i,j+\frac{1}{2}}),  & \mbox{if} ~~2\le j \le n_t-1,
		\\
		\frac{1}{2}\rho^{-1}_{i,\frac{3}{2}}, & \mbox{if} ~~j=1,
		\\
		\frac{1}{2}\rho^{-1}_{i,n_t-\frac{1}{2}}, & \mbox{if} ~~j=n_t,
		\end{cases}
	\]
	\[
		a_{i,j}=
		\begin{cases}
		\frac{1}{2}(\rho^{0}_i)^{-1} & \mbox{if}~~j=1,\\
		\frac{1}{2}(\rho^{1}_i)^{-1} & \mbox{if}~~j=n_t,\\
		0 & \mbox{otherwise}.
		\end{cases}
	\]
We remark that it is important to first square then average, and first invert then average, to guarantee stability \cite{Asc08,HabHor15}.

\subsection{Sequential quadratic programming (SQP)}

Following the above discretization scheme, we obtain the discrete convex optimization problem
	\begin{subequations}\label{eq:SQPformulation}
	\begin{eqnarray}
		\min && f(p, \rho, u) =\left\langle A_1 (p^*p), A_2 (\rho^{-1})+a\right\rangle h_xh_t+
		\left\langle u^*u, A_2 (\rho^{-1})+a\right\rangle\gamma h_xh_t,\\
		{\rm s.t.} && D_1 p+D_2 \rho+D_3 u=b.
	\end{eqnarray}
	\end{subequations}
The Lagrangian of this problem is
	\[
		\cL(p, \rho, u)= f(p, \rho, u)/(h_xh_t)+\left\langle \lambda, D_1 p+D_2 \rho+D_3 u-b\right\rangle.
	\]
The KKT condition \cite{ByrCurNoc08,NocWri06}
	\begin{subequations}
	\begin{eqnarray}
	\nabla_p \cL &=& D_1^* \lambda+2p\circ A_1^*(A_2(\rho^{-1})+a)=0
	\\
	\nabla_\rho \cL &=& D_2^* \lambda-\rho^{-1}\circ A_2^*A_1(p^*p)\circ \rho^{-1}
	-\gamma\rho^{-1}\circ A_2^*(u^*u)\circ \rho^{-1}=0
	\\
	\nabla_u \cL &=& D_3^*\lambda+2\gamma u\circ (A_2(\rho^{-1})+a)=0
	\\
	\nabla_\lambda\cL &=& D_1 p+D_2 \rho+D_3 u-b=0
	\end{eqnarray}
	\end{subequations}
follow, with $\circ$ denoting block-wise multiplication.

Let $w=(p,\rho, u)$, $D=(D_1, D_2, D_3)$, then at each SQP iteration we solve the system
	\begin{equation}\label{eq:SQPiter}
		\left (\begin{matrix}\hat A & D^* \\ D &0\end{matrix}\right)
		\left(\begin{array}{c}\delta w\\\delta \lambda\end{array}\right)=
		-\left(\begin{array}{c}\nabla_w\cL \\ \nabla_\lambda \cL\end{array}\right),
	\end{equation}
and update $w,\lambda$ using line search. In principle, Problem~\ref{eq:SQPformulation} can be solved using Newton's method.
However, the mixed terms introduce off-diagonal elements in the Hessian, which makes it forbidden for large problems.
We adopt an inexact SQP method \cite{ByrCurNoc08}. The matrix $\hat A$ is an approximation of the Hessian of the objective function
	\[
		\hat A=
		\left(\begin{matrix}
		2 {\rm Bdiag} (A_1^* (A_2(\rho^{-1})+a)) & 0 & 0
		\\
		0 & {\rm Bdiag} (g(p,\rho,u)) & 0
		\\
		0 & 0 & 2 \gamma{\rm Bdiag} (A_2(\rho^{-1})+a)
		\end{matrix}\right).
	\]
Here ${\rm Bdiag}$ denotes block diagonal operator. More specifically,
	\[
		{\rm Bdiag}(T_1, T_2,\cdots, T_k)=
		\left[\begin{matrix}
		T_1 & 0 & \cdots & 0\\
		0 & T_2 & \cdots & 0\\
		\vdots & \vdots & \ddots & \vdots\\
		0 & 0 & \cdots & T_k
		\end{matrix}\right]
	\]
for linear operators $T_1, T_2, \cdots, T_k$.	
The operator $g(p, \rho, u)$ is the Hessian of $f$ over $\rho$ with $g_{i,j+\frac{1}{2}}$ being the map
	\begin{eqnarray*}
		g_{i,j+\frac{1}{2}}(X)&=&\rho^{-1}_{i,j+\frac{1}{2}}(A_2^*A_1(p^*p))_{i,j+\frac{1}{2}}
		\rho^{-1}_{i,j+\frac{1}{2}}X\rho^{-1}_{i,j+\frac{1}{2}}+\rho^{-1}_{i,j+\frac{1}{2}}X
		\rho^{-1}_{i,j+\frac{1}{2}}(A_2^*A_1(p^*p))_{i,j+\frac{1}{2}}\rho^{-1}_{i,j+\frac{1}{2}}
		\\
		&&+ \gamma\rho^{-1}_{i,j+\frac{1}{2}}(A_2^*(u^*u))_{i,j+\frac{1}{2}}
		\rho^{-1}_{i,j+\frac{1}{2}}X\rho^{-1}_{i,j+\frac{1}{2}}+\gamma\rho^{-1}_{i,j+\frac{1}{2}}X
		\rho^{-1}_{i,j+\frac{1}{2}}(A_2^*(u^*u))_{i,j+\frac{1}{2}}\rho^{-1}_{i,j+\frac{1}{2}}.
	\end{eqnarray*}
	
In each step we solve the linear system \eqref{eq:SQPiter} in an inexact manner. There are many methods to achieve this. In our approach, we apply the
Schur complement and solve the reduced system
	\[
		D\hat A^{-1} D^* \delta \lambda = \nabla_\lambda \cL- D\hat A^{-1} \nabla_w \cL
	\]
using preconditioned conjugated gradients method with incomplete Cholesky factorization \cite{Ker78} as a preconditioner. The update for $w$ is then given by
	\[
		\delta w = -\hat A^{-1} (D^*\delta\lambda+\nabla_w\cL).
	\]
\begin{remark}
In our numerical implementation, we take advantage of the structure of $\rho$ being symmetric, and only save the upper triangular part of it.
This is beneficial in terms of both memory and  speed.
\end{remark}

\subsection{2D and 3D cases} \label{sec:2D}

In this section we sketch what happens in higher dimensions, namely 2D and 3D.

We begin with the 2D case. Accordingly, we have the discrete convex optimization problem
	\begin{eqnarray*}
		\min && f(p, \rho, u) =\left\langle A_{1x} (p_x^*p_x)+A_{1y}(p_y^*p_y), A_2 (\rho^{-1})+
		a\right\rangle h_xh_yh_t+
		\left\langle u^*u, A_2 (\rho^{-1})+a\right\rangle\gamma h_xh_yh_t\\
		{\rm s.t.} && D_{1x} p_x +D_{1y} p_y+D_2 \rho+D_3 u=b.
	\end{eqnarray*}
The Lagrangian of this problem is
	\[
		\cL(p, \rho, u)= f(p, \rho, u)/(h_xh_yh_t)+\left\langle \lambda,
		D_{1x} p_x+D_{1y} p_y+D_2 \rho+D_3 u-b\right\rangle.
	\]
In the above,
	\[
		a_{i,j,k}=
		\begin{cases}
		\frac{1}{2}(\rho^{0}_{i,j})^{-1} & \mbox{if}~~k=1,\\
		\frac{1}{2}(\rho^{1}_{i,j})^{-1} & \mbox{if}~~k=n_t,\\
		0 & \mbox{otherwise}.
		\end{cases}
	\]
and
	\[
		b_{i,j,k}=
		\begin{cases}
		\rho^0_{i,j}/h_t & \mbox{if}~~k=1,\\
		-\rho^1_{i,j}/h_t & \mbox{if}~~k=n_t,\\
		0 & \mbox{otherwise}.
		\end{cases}
	\]
It follows that the KKT conditions are
	\begin{subequations}
	\begin{eqnarray}
	\nabla_{p_x} \cL &=& D_{1x}^* \lambda+2p_x\circ A_{1x}^*(A_2(\rho^{-1})+a)=0
	\\
	\nabla_{p_y} \cL &=& D_{1y}^* \lambda+2p_y\circ A_{1y}^*(A_2(\rho^{-1})+a)=0
	\\
	\nabla_\rho \cL &=& D_2^* \lambda-\rho^{-1}\circ A_2^*(A_{1x}(p_x^*p_x)+A_{1y}(p_y^*p_y))
	\circ \rho^{-1}
	-\gamma\rho^{-1}\circ A_2^*(u^*u)\circ \rho^{-1}=0
	\\
	\nabla_u \cL &=& D_3^*\lambda+2\gamma u\circ (A_2(\rho^{-1})+a)=0
	\\
	\nabla_\lambda\cL &=& D_1 p+D_2 \rho+D_3 u-b=0,
	\end{eqnarray}
	\end{subequations}
with $\circ$ denoting block-wise multiplication as before.

Let $w=(p_x,p_y,\rho, u).$ Then at each SQP iteration, we solve the system
	\begin{equation}
		\left (\begin{matrix}\hat A & D^* \\ D &0\end{matrix}\right)
		\left(\begin{array}{c}\delta w\\\delta \lambda\end{array}\right)=
		-\left(\begin{array}{c}\nabla_w\cL \\ \nabla_\lambda \cL\end{array}\right),
	\end{equation}
where $D=(D_{1x},D_{1y}, D_2, D_3)$. The matrix $\hat A$ is an approximation of the Hessian of the objective function
	\[
		\hat A=
		\left(\begin{matrix}
		2 {\rm Bdiag} (A_{1x}^* (A_2(\rho^{-1})+a)) & 0 & 0 & 0
		\\
		0 & 2 {\rm Bdiag} (A_{1y}^* (A_2(\rho^{-1})+a)) & 0 & 0
		\\
		0 & 0 & {\rm Bdiag} (g(p,\rho,u)) & 0
		\\
		0 & 0 & 0 & 2 \gamma{\rm Bdiag} (A_2(\rho^{-1})+a)
		\end{matrix}\right)
	\]
The operator $g(p, \rho, u)$ is the Hessian of $f$ over $\rho$ with $g_{i,j,k+\frac{1}{2}}$ being the map
	\begin{eqnarray*}
		g_{i,j,k+\frac{1}{2}}(X)&=&\rho^{-1}_{i,j,k+\frac{1}{2}}(A_2^*(A_{1x}(p_x^*p_x)+
		A_{1y}(p_y^*p_y)+\gamma u^*u))_{i,j,k+\frac{1}{2}}
		\rho^{-1}_{i,j,k+\frac{1}{2}}X\rho^{-1}_{i,j,k+\frac{1}{2}}
		\\&&+\rho^{-1}_{i,j,k+\frac{1}{2}}X
		\rho^{-1}_{i,j,k+\frac{1}{2}}(A_2^*(A_{1x}(p_x^*p_x)+
		A_{1y}(p_y^*p_y)+\gamma u^*u))_{i,j,k+\frac{1}{2}}\rho^{-1}_{i,j,k+\frac{1}{2}}.
	\end{eqnarray*}

The 3D case is quite similar. Now, we have the discrete convex optimization problem
	\begin{eqnarray*}
		\min && f(p, \rho, u) =\left\langle A_{1x} (p_x^*p_x)+A_{1y}(p_y^*p_y)+A_{1z}(p_z^*p_z),
		 A_2 (\rho^{-1})+
		a\right\rangle h_xh_yh_zh_t\\&&\hspace{1.9cm}+
		\left\langle u^*u, A_2 (\rho^{-1})+a\right\rangle\gamma h_xh_yh_zh_t\\
		{\rm s.t.} && D_{1x} p_x +D_{1y} p_y+D_{1z}p_z+D_2 \rho+D_3 u=b.
	\end{eqnarray*}
The Lagrangian of this problem is
	\[
		\cL(p, \rho, u)= f(p, \rho, u)/(h_xh_yh_zh_t)+\left\langle \lambda,
		D_{1x} p_x+D_{1y} p_y+D_{1z}p_z+D_2 \rho+D_3 u-b\right\rangle.
	\]
In the above,
	\[
		a_{i,j,k,\ell}=
		\begin{cases}
		\frac{1}{2}(\rho^{0}_{i,j,k})^{-1} & \mbox{if}~~\ell=1,\\
		\frac{1}{2}(\rho^{1}_{i,j,k})^{-1} & \mbox{if}~~\ell=n_t,\\
		0 & \mbox{otherwise}.
		\end{cases}
	\]
and
	\[
		b_{i,j,k,\ell}=
		\begin{cases}
		\rho^0_{i,j,k}/h_t & \mbox{if}~~\ell=1,\\
		-\rho^1_{i,j,k}/h_t & \mbox{if}~~\ell=n_t,\\
		0 & \mbox{otherwise}.
		\end{cases}
	\]
It follows that the KKT conditions now are
	\begin{eqnarray*}
	\nabla_{p_x} \cL &=& D_{1x}^* \lambda+2p_x\circ A_{1x}^*(A_2(\rho^{-1})+a)=0
	\\
	\nabla_{p_y} \cL &=& D_{1y}^* \lambda+2p_y\circ A_{1y}^*(A_2(\rho^{-1})+a)=0
	\\
	\nabla_{p_z} \cL &=& D_{1z}^* \lambda+2p_z\circ A_{1z}^*(A_2(\rho^{-1})+a)=0
	\\
	\nabla_\rho \cL &=& D_2^* \lambda-\rho^{-1}\circ A_2^*(A_{1x}(p_x^*p_x)+A_{1y}(p_y^*p_y)
	+A_{1z}(p_z^*p_z))
	\circ \rho^{-1}
	-\gamma\rho^{-1}\circ A_2^*(u^*u)\circ \rho^{-1}=0
	\\
	\nabla_u \cL &=& D_3^*\lambda+2\gamma u\circ (A_2(\rho^{-1})+a)=0
	\\
	\nabla_\lambda\cL &=& D_1 p+D_2 \rho+D_3 u-b=0,
	\end{eqnarray*}
with $\circ$ the block-wise multiplication as earlier.

Let $w=(p_x,p_y,p_z,\rho, u)$, then at each SQP iteration we solve the system
	\begin{equation}
		\left (\begin{matrix}\hat A & D^* \\ D &0\end{matrix}\right)
		\left(\begin{array}{c}\delta w\\\delta \lambda\end{array}\right)=
		-\left(\begin{array}{c}\nabla_w\cL \\ \nabla_\lambda \cL\end{array}\right),
	\end{equation}
where $D=(D_{1x},D_{1y},D_{1z}, D_2, D_3)$. The matrix $\hat A$ is an approximation of the Hessian of the objective function
	{\scriptsize
	\[
		\left(\begin{matrix}
		2 {\rm Bdiag} (A_{1x}^* (A_2(\rho^{-1})+a)) & 0 & 0 & 0 & 0
		\\
		0 & 2 {\rm Bdiag} (A_{1y}^* (A_2(\rho^{-1})+a)) & 0 & 0 & 0
		\\
		0 & 0 & 2 {\rm Bdiag} (A_{1z}^* (A_2(\rho^{-1})+a)) & 0 & 0
		\\
		0 & 0 & 0 & {\rm Bdiag} (g(p,\rho,u)) & 0
		\\
		0 & 0 & 0 & 0 & 2 \gamma{\rm Bdiag} (A_2(\rho^{-1})+a)
		\end{matrix}\right)
	\]
	}
The operator $g(p, \rho, u)$ is the Hessian of $f$ over $\rho$ with $g_{i,j,k,\ell+\frac{1}{2}}$ being the map
	\begin{eqnarray*}
		g_{i,j,k,\ell+\frac{1}{2}}(X)&=&\rho^{-1}_{i,j,k,\ell+\frac{1}{2}}(A_2^*(A_{1x}(p_x^*p_x)+
		A_{1y}(p_y^*p_y)+A_{1z}(p_z^*p_z)+\gamma u^*u))_{i,j,k,\ell+\frac{1}{2}}
		\rho^{-1}_{i,j,k,\ell+\frac{1}{2}}X\rho^{-1}_{i,j,k,\ell+\frac{1}{2}}
		\\&&+\rho^{-1}_{i,j,k,\ell+\frac{1}{2}}X
		\rho^{-1}_{i,j,k,\ell+\frac{1}{2}}(A_2^*(A_{1x}(p_x^*p_x)+A_{1y}(p_y^*p_y)+A_{1z}(p_z^*p_z)
		+\gamma u^*u))_{i,j,k,\ell+\frac{1}{2}}\rho^{-1}_{i,j,k,\ell+\frac{1}{2}}
	\end{eqnarray*}

\section{Vector-valued optimal mass transport} \label{sec:vector}

Next we move to vector-valued optimal transport, which was proposed recently in \cite{CheGeoTan17}. We briefly review the setup in this section, and
refer the reader to \cite{CheGeoTan17} for details.

\subsection{Gradients on graphs}

We consider a connected, positively weighted, undirected graph $\cF =(\cV, \cE,\cW)$ with $n$ nodes labeled as $i$, with $1\le i \le n$, and $N$ edges.
We have that $\Delta_\cF =-\mD W\mD^T$ where $\Delta_\cF, \mD, W=\diag \{\bw_1, \cdots, \bw_N\}$ are the graph Laplacian, incidence, and weight
matrices, respectively. One can define the Laplacian in terms of a graph gradient and divergence as
\[  \Delta_\cF = -\nabla_\cF^*\nabla_\cF,
\]
where
	\[
		\nabla_\cF: \mR^n \rightarrow \mR^N, ~ x \mapsto W^{1/2}\mD^T x
	\]
denotes the gradient operator and
	\[
		\nabla_\cF^*:  \mR^N \rightarrow \mR^n, ~y \mapsto \mD W^{1/2} y
	\]
denotes its dual.

\subsection{Vector-valued optimal mass transport}

We begin by considering a \textbf{\emph{vector-valued density}} $\rho$ on $\mR^m$, i.e., a map from $E\subset \mR^m$ to $\mR_+^{n}$ such that
	\[
		\sum_{i=1}^n\int_E \rho_i(x)dx=1.
	\]
Here the convex compact set $E \subset \mR^m$ is a domain where the densities are defined, typically the unit $n$-dimensional cube.
To avoid proliferation of symbols, we denote the set of all vector-valued densities and its interior again by $\cD$ and $\cD_+$, respectively.
We refer to the entries of $\rho$ as representing density or mass of individual species/particles that can mutate between one another while maintaining total mass.
Mass transfer may only be permissible between specific types of particles. Thus, allowable transfer can be modeled by the existence of a corresponding edge in a graph $\cF=(\cV,\cE,\cW)$ whose vertices in $\cV$ correspond to those individual species, see \cite{CheGeoTan17}. The edge weights in $\cW$ can quantify cost, rate, or likelihood of transfer.

Following the arguments in \cite{CheGeoTan17}, this leads to the following (symmetric) Wasssertein 2-metric on vector-valued distributions:
Given two given marginals $\rho^0, \rho^1 \in \cD_+$ the (square) of the Wasserstein distance is given by:
	 \begin{subequations}\label{eq:vecomtcvx}
 	\begin{eqnarray}\label{eq:vecomtcvx1}
	&&\min_{\rho, p, u} \int_0^1\int_{E}
	\left \{p^T\diag(\rho)^{-1}p+\gamma u^T[\diag(\mD_2^T\rho)^{-1}+
	\diag(\mD_1^T\rho)^{-1}]u\right\}dxdt\\
	&&\frac{\partial \rho}{\partial t}+\nabla_x\cdot p-\nabla_\cF^*\, u=0, \label{eq:vecomtcvx2}\\
	&& \rho(0,\cdot)=\rho^0, ~~\rho(1,\cdot)=\rho^1  \label{eq:vecomtcvx3}.
	\end{eqnarray}
	\end{subequations}
Here $u$ is the ``flux'' on graphs, $p=[p_1, \cdots, p_n]^T$ is the ``momentum'' (mass times velocity vector field), the matrix $\mD_1$ is the portion
of the incidence matrix $\mD$ containing 1's (sources), and $\mD_2 = \mD_1-\mD$ (sinks).
In what follows, we describe an algorithm for the numerical implementation of this convex optimization problem.

\section{Discretization and algorithm: vector-valued case}

As in the matrix-valued cases, for simplicity of exposition, we consider the discretization in 1D case, and describe the 2D case in Section~\ref{sec:2D_vector} below.
Thus, we take $E=[0,1]$, and as before our technique extends almost verbatim to the higher dimensional setting; see Section~\ref{sec:2D_vector}.
We should note that the algorithm presented here in the vector-valued case is very similar to the matrix optimal transport just described in the preceding sections.

We discretize the space-time domain $[0,\,1]\times[0,\,1]$ into $n_x\times n_t$ rectangular cells. Denote $\Omega_{ij}, 1\le i\le n_x, 1\le j\le n_t$ as the $(i,j)$ box.
We use staggered grid to discretize $p$ and $\rho$. {\bf \em The variable $u$ is, however, valued at the centers of the cells $\{\Omega_{ij}\}$}. More specifically,
	\[
		p=(p_{i+\frac{1}{2},j}), ~0\le i\le n_x, ~1\le j\le n_t
	\]
	\[
		\rho=(\rho_{i,j+\frac{1}{2}}),~1\le i\le n_x,~0\le j\le n_t
	\]
	\[
		u=(u_{i,j}),~1\le i \le n_x, ~1\le j\le n_t.
	\]
Note that the boundary values are
	\[
		p_{\frac{1}{2},j}=0, ~p_{n_x+\frac{1}{2},j}=0, ~1\le j\le n_t
	\]
and
	\[
		\rho_{i,\frac{1}{2}}=\rho^0_i, ~\rho_{i,n_t+\frac{1}{2}}=\rho^1_i, ~1\le i\le n_x.
	\]
We exclude the boundary values from the variables and denote
	\[
		p=(p_{i+\frac{1}{2},j}), ~1\le i\le n_x-1, ~1\le j\le n_t
	\]
	\[
		\rho=(\rho_{i,j+\frac{1}{2}}),~1\le i\le n_x,~1\le j\le n_t-1.
	\]
	
\subsection{Continuity equation}

We use the preceding discretizing scheme, together with the boundary conditions to rewrite the continuity equation \eqref{eq:vecomtcvx2} as
	\begin{equation}\label{eq:vecdiscretecontinuity}
	D_1 p+D_2 \rho+D_3 u=b.
	\end{equation}
Here the linear operators $D_1, D_2, D_3$ are defined as
	\[
		(D_1 p)_{i,j}=
		\begin{cases}
		(p_{i+\frac{1}{2},j}
		-p_{i-\frac{1}{2},j})/h_x,  & \mbox{if} ~~2\le i \le n_x-1,
		\\
		p_{\frac{3}{2},j}/h_x, & \mbox{if} ~~i=1,
		\\
		-p_{n_x-\frac{1}{2},j}/h_x, & \mbox{if} ~~i=n_x,
		\end{cases}
	\]
	\[
		(D_2 \rho)_{i,j}=
		\begin{cases}
		(\rho_{i,j+\frac{1}{2}}-\rho_{i,j-\frac{1}{2}})/h_t,  & \mbox{if} ~~2\le j \le n_t-1,
		\\
		 \rho_{i,\frac{3}{2}}/h_t, & \mbox{if} ~~j=1,
		\\
		- \rho_{i,n_t-\frac{1}{2}}/h_t, & \mbox{if} ~~j=n_t,
		\end{cases}
	\]
	\[
		(D_3 u)_{i,j}=-\nabla_\cF^* u_{i,j}, ~1\le i \le n_x, ~1\le j\le n_t.
	\]
The parameter $b$ carries the information of the boundary values $\rho^0$ and $\rho^1$. More specifically,
	\[
		b_{i,j}=
		\begin{cases}
		\rho^0_i/h_t & \mbox{if}~~j=1,\\
		-\rho^1_i/h_t & \mbox{if}~~j=n_t,\\
		0 & \mbox{otherwise}.
		\end{cases}
	\]

\subsection{Discretization of the cost function}

Let $A_1$ be the averaging operator over the spatial domain and $A_2$ be the averaging operator over the time domain (as before one needs to be careful about the boundaries).
Then the cost function  \eqref{eq:vecomtcvx1} may be approximated by
	\begin{equation}
		\left\langle A_1 (p^2), A_2 (1/\rho)+a\right\rangle h_xh_t+
		\left\langle u^2, A_2 (1/(\mD_2^T\rho)+1/(\mD_1^T\rho))+c\right\rangle\gamma h_xh_t,
	\end{equation}
where $a\ge 0$ depends only on the boundary values $\rho^0$ and $\rho^1$. The inverse operator and multiplication operators are applied block-wise.
The expressions for $A_1, A_2, a$ are
	\[
		(A_1 (p^2))_{i,j}=
		\begin{cases}
		\frac{1}{2}(p^2_{i-\frac{1}{2},j}
		+p^2_{i+\frac{1}{2},j}),  & \mbox{if} ~~2\le i \le n_x-1,
		\\
		\frac{1}{2}p^2_{\frac{3}{2},j}, & \mbox{if} ~~i=1,
		\\
		\frac{1}{2}p^2_{n_x-\frac{1}{2},j}, & \mbox{if} ~~i=n_x,
		\end{cases}
	\]
	\[
		(A_2 (1/\rho))_{i,j}=
		\begin{cases}
		\frac{1}{2}(1/\rho_{i,j-\frac{1}{2}}
		+1/\rho_{i,j+\frac{1}{2}}),  & \mbox{if} ~~2\le j \le n_t-1,
		\\
		1/\rho_{i,\frac{3}{2}}/2, & \mbox{if} ~~j=1,
		\\
		1/\rho_{i,n_t-\frac{1}{2}}/2, & \mbox{if} ~~j=n_t,
		\end{cases}
	\]
	\[
		a_{i,j}=
		\begin{cases}
		1/\rho^{0}_i/2& \mbox{if}~~j=1,\\
		1/\rho^{1}_i/2 & \mbox{if}~~j=n_t,\\
		0 & \mbox{otherwise},
		\end{cases}
	\]
	\[
		c_{i,j}=
		\begin{cases}
		1/\mD_2^T\rho^{0}_i/2+1/\mD_1^T\rho^{0}_i/2& \mbox{if}~~j=1,\\
		1/\mD_2^T\rho^{1}_i/2+1/\mD_1^T\rho^{1}_i/2 & \mbox{if}~~j=n_t,\\
		0 & \mbox{otherwise}.
		\end{cases}
	\]

\subsection{Sequential quadratic programming (SQP)}

From the above discussion, we obtain the discrete convex optimization problem
	\begin{subequations}
	\begin{eqnarray}
		\min && f(p, \rho, u) =\left\langle A_1 (p^2), A_2 (1/\rho)+a\right\rangle h_xh_t+
		\left\langle u^2, A_2 (1/(\mD_2^T\rho)+1/(\mD_1^T\rho))+c\right\rangle\gamma h_xh_t\\
		{\rm s.t.} && D_1 p+D_2 \rho+D_3 u=b.
	\end{eqnarray}
	\end{subequations}
The Lagrangian of this problem is
	\[
		\cL(p, \rho, u)= f(p, \rho, u)/(h_xh_t)+\left\langle \lambda, D_1 p+D_2 \rho+D_3 u-b\right\rangle.
	\]
It follows that the KKT conditions are given by
	\begin{subequations}
	\begin{eqnarray}
	\nabla_p \cL &=& D_1^T \lambda+2p\circ A_1^T(A_2(1/\rho)+a)=0
	\\
	\nabla_\rho \cL &=& D_2^T \lambda-A_2^TA_1(p^2)/ \rho^{2}
	-\gamma \mD_2(A_2^T(u^2)/(\mD_2^T\rho)^2)-\gamma \mD_1(A_2^T(u^2)/(\mD_1^T\rho)^2)=0
	\\
	\nabla_u \cL &=& D_3^T\lambda+2\gamma u\circ (A_2(1/(\mD_2^T\rho)+1/(\mD_1^T\rho))+c)=0
	\\
	\nabla_\lambda\cL &=& D_1 p+D_2 \rho+D_3 u-b=0,
	\end{eqnarray}
	\end{subequations}
with $\circ$ denoting block-wise multiplication.

Let $w=(p,\rho, u).$ Then at each SQP iteration, we solve the system
	\begin{equation}
		\left (\begin{matrix}\hat A & D^T \\ D &0\end{matrix}\right)
		\left(\begin{array}{c}\delta w\\\delta \lambda\end{array}\right)=
		-\left(\begin{array}{c}\nabla_w\cL \\ \nabla_\lambda \cL\end{array}\right),
	\end{equation}
where $D=(D_1, D_2, D_3)$. Again, the matrix $\hat A$ is an approximation of the Hessian of the objective function
	\[
		\hat A=
		\left(\begin{matrix}
		2 {\rm diag} (A_1^T (A_2(1/\rho)+a)) & 0 & 0
		\\
		0 & {\rm diag} (g(p,\rho,u)) & 0
		\\
		0 & 0 & 2 \gamma{\rm diag} (A_2(1/(\mD_2^T\rho)+1/(\mD_1^T\rho))+c)
		\end{matrix}\right).
	\]
The operator $g(p, \rho, u)$ is the Hessian of $f$ over $\rho$ with $g_{i,j+\frac{1}{2}}$ being the map
	\begin{eqnarray*}
		g_{i,j+\frac{1}{2}}(X)&=&2(A_2^TA_1(p^2))_{i,j+\frac{1}{2}}
		/\rho^{3}_{i,j+\frac{1}{2}}X\\
		&+& 2\gamma\mD_2[(A_2^T(u^2))_{i,j+\frac{1}{2}}
		/(\mD_2^T\rho)^{3}_{i,j+\frac{1}{2}}\mD_2^TX]+2\gamma\mD_1[(A_2^T(u^2))_{i,j+\frac{1}{2}}
		/(\mD_1^T\rho)^{3}_{i,j+\frac{1}{2}}\mD_1^TX].
	\end{eqnarray*}

\subsection{2D case} \label{sec:2D_vector}

We concretely work out the 2D case in this section. The higher dimensional cases are very similar, but naturally involve additional indices.
We have the discrete convex optimization problem
	\begin{eqnarray*}
		\min && f(p, \rho, u) =\left\langle A_{1x} (p_x^2)+A_{1y}(p_y^2), A_2 (1/\rho)+
		a\right\rangle h_xh_yh_t+
		\left\langle u^2, A_2 (1/(\mD_2^T\rho)+1/(\mD_1^T\rho))+c\right\rangle\gamma h_xh_yh_t\\
		{\rm s.t.} && D_{1x} p_x +D_{1y} p_y+D_2 \rho+D_3 u=b.
	\end{eqnarray*}
The Lagrangian of this problem is
	\[
		\cL(p, \rho, u)= f(p, \rho, u)/(h_xh_yh_t)+\left\langle \lambda,
		D_{1x} p_x+D_{1y} p_y+D_2 \rho+D_3 u-b\right\rangle.
	\]
In the above,
	\[
		a_{i,j,k}=
		\begin{cases}
		1/\rho^{0}_{i,j}/2& \mbox{if}~~k=1,\\
		1/\rho^{1}_{i,j}/2 & \mbox{if}~~k=n_t,\\
		0 & \mbox{otherwise},
		\end{cases}
	\]
and
	\[
		b_{i,j,k}=
		\begin{cases}
		\rho^0_{i,j}/h_t & \mbox{if}~~k=1,\\
		-\rho^1_{i,j}/h_t & \mbox{if}~~k=n_t,\\
		0 & \mbox{otherwise},
		\end{cases}
	\]
	\[
		c_{i,j,k}=
		\begin{cases}
		1/\mD_2^T\rho^{0}_{i,j}/2+1/\mD_1^T\rho^{0}_{i,j}/2& \mbox{if}~~k=1,\\
		1/\mD_2^T\rho^{1}_{i,j}/2+1/\mD_1^T\rho^{1}_{i,j}/2 & \mbox{if}~~k=n_t,\\
		0 & \mbox{otherwise}.
		\end{cases}
	\]
The KKT conditions now are
\begin{eqnarray*}
	\nabla_{p_x} \cL &=& D_{1x}^T \lambda+2p_x\circ A_{1x}^T(A_2(1/\rho)+a)=0
	\\
	\nabla_{p_y} \cL &=& D_{1y}^T \lambda+2p_y\circ A_{1y}^T(A_2(1/\rho)+a)=0
	\\
	\nabla_\rho \cL &=& D_2^T \lambda-A_2^T(A_{1x}(p_x^2)+A_{1y}(p_y^2))/\rho^2
	-\gamma \mD_2(A_2^T(u^2)/(\mD_2^T\rho)^2)-\gamma \mD_1(A_2^T(u^2)/(\mD_1^T\rho)^2)=0
	\\
	\nabla_u \cL &=& D_3^T\lambda+2\gamma u\circ (A_2(1/(\mD_2^T\rho)+1/(\mD_1^T\rho))+c)=0
	\\
	\nabla_\lambda\cL &=& D_1 p+D_2 \rho+D_3 u-b=0,
	\end{eqnarray*}
with $\circ$ denoting block-wise multiplication.

Let $w=(p_x,p_y,\rho, u)$, then at each SQP iteration we solve the system
	\begin{equation}
		\left (\begin{matrix}\hat A & D^* \\ D &0\end{matrix}\right)
		\left(\begin{array}{c}\delta w\\\delta \lambda\end{array}\right)=
		-\left(\begin{array}{c}\nabla_w\cL \\ \nabla_\lambda \cL\end{array}\right),
	\end{equation}
where $D=(D_{1x},D_{1y}, D_2, D_3)$. The matrix $\hat A$ is an approximation of the Hessian of the objective function
{\footnotesize
	\[
		\hspace{-1cm}\left(\begin{matrix}
		2 {\rm diag} (A_{1x}^T (A_2(1/\rho)+a)) & 0 & 0 & 0
		\\
		0 & 2 {\rm diag} (A_{1y}^T (A_2(1/\rho)+a)) & 0 & 0
		\\
		0 & 0 & {\rm diag} (g(p,\rho,u)) & 0
		\\
		0 & 0 & 0 & 2 \gamma{\rm diag} (A_2(1/(\mD_2^T\rho)+1/(\mD_1^T\rho))+c)
		\end{matrix}\right).
	\]
	}
The operator $g(p, \rho, u)$ is the Hessian of $f$ over $\rho$ with $g_{i,j,k+\frac{1}{2}}$ being the map
	\begin{eqnarray*}
		g_{i,j+\frac{1}{2}}(X)&=&2(A_2^T(A_{1x}(p_x^2)+A_{1y}(p_y^2)))_{i,j+\frac{1}{2}}
		/\rho^{3}_{i,j+\frac{1}{2}}X\\
		&+& 2\gamma\mD_2[(A_2^T(u^2))_{i,j+\frac{1}{2}}
		/(\mD_2^T\rho)^{3}_{i,j+\frac{1}{2}}\mD_2^TX]+2\gamma\mD_1[(A_2^T(u^2))_{i,j+\frac{1}{2}}
		/(\mD_1^T\rho)^{3}_{i,j+\frac{1}{2}}\mD_1^TX].
	\end{eqnarray*}

\section{Numerical experiments} \label{sec:numerical}

Several examples are provided in this section to illustrate the effectiveness of our algorithms.
For matrix-valued densities, we present examples in both 2D and 3D settings. In contrast, only 2D examples are studied for vector-valued densities.

\subsection{Matrix case}
One motivation for matrix-valued optimal mass transport comes from diffusion tensor imaging (DTI). This is a widely used technique in magnetic resonance imaging.
In diffusion images, the information at each pixel is captured in a ellipsoid, i.e., a $3 \times 3$ positive definite matrix,
in lieu of a nonnegative number. The ellipsoids describe useful information such as the orientations of the brain fibers.

\begin{figure}[h]
\centering
\subfloat[$\rho^0$]{\includegraphics[width=0.40\textwidth]{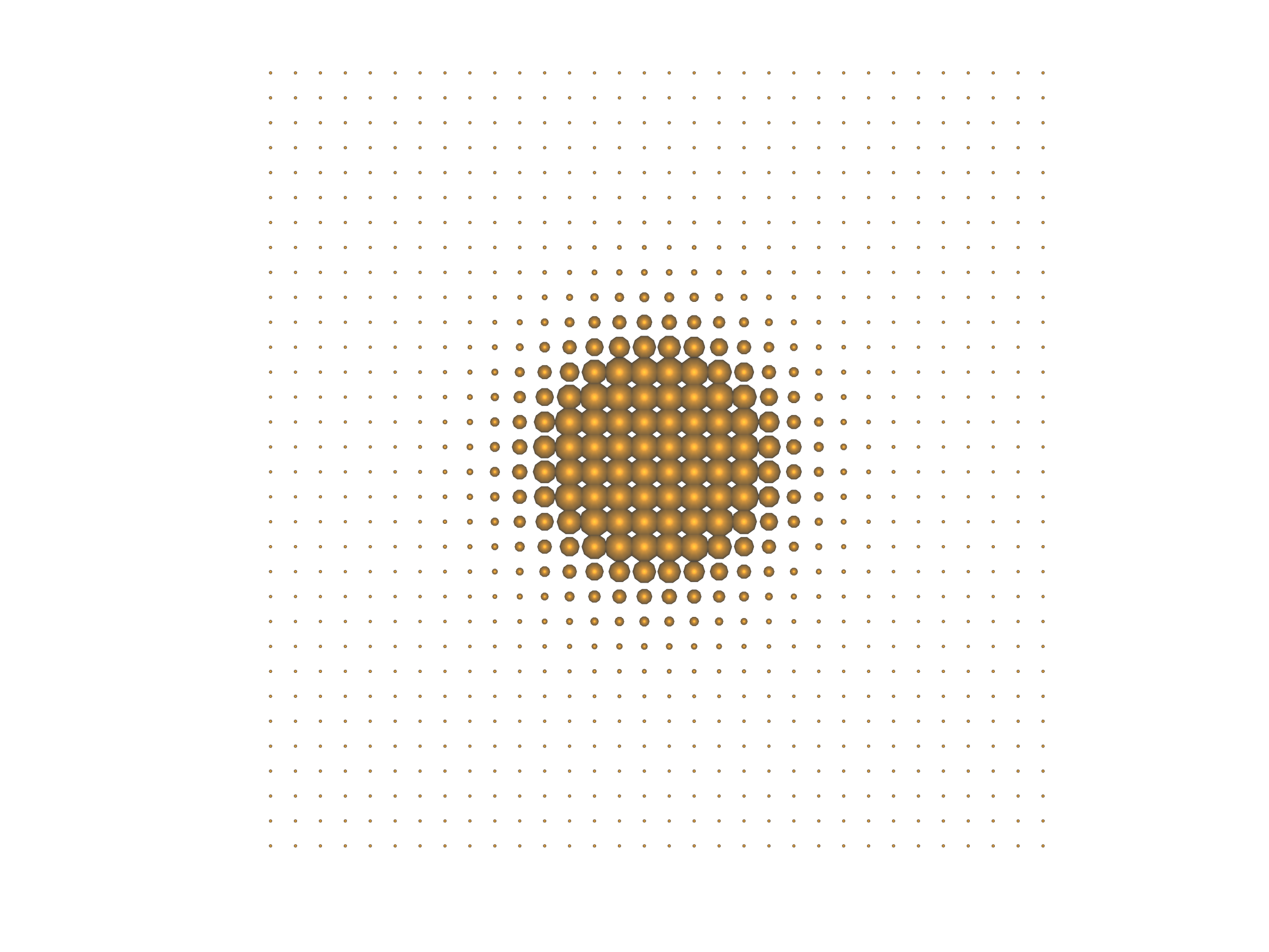}}
\subfloat[$\rho^1$]{\includegraphics[width=0.40\textwidth]{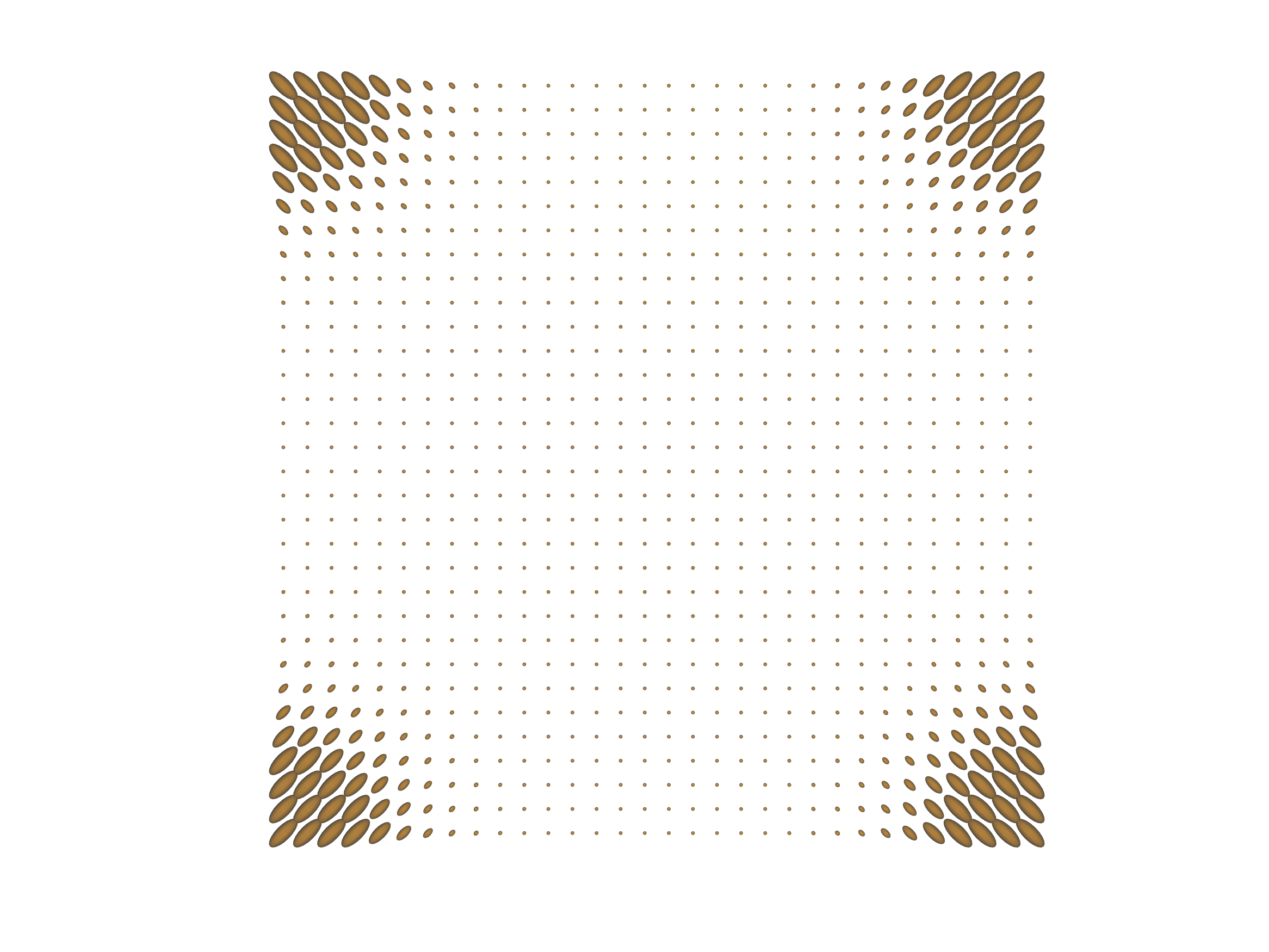}}
 \caption{Marginal distributions}
 \label{fig:ex1marginals}
\end{figure}

\begin{figure}[h]
\centering
\subfloat[$t=0.1$]{\includegraphics[width=0.34\textwidth]{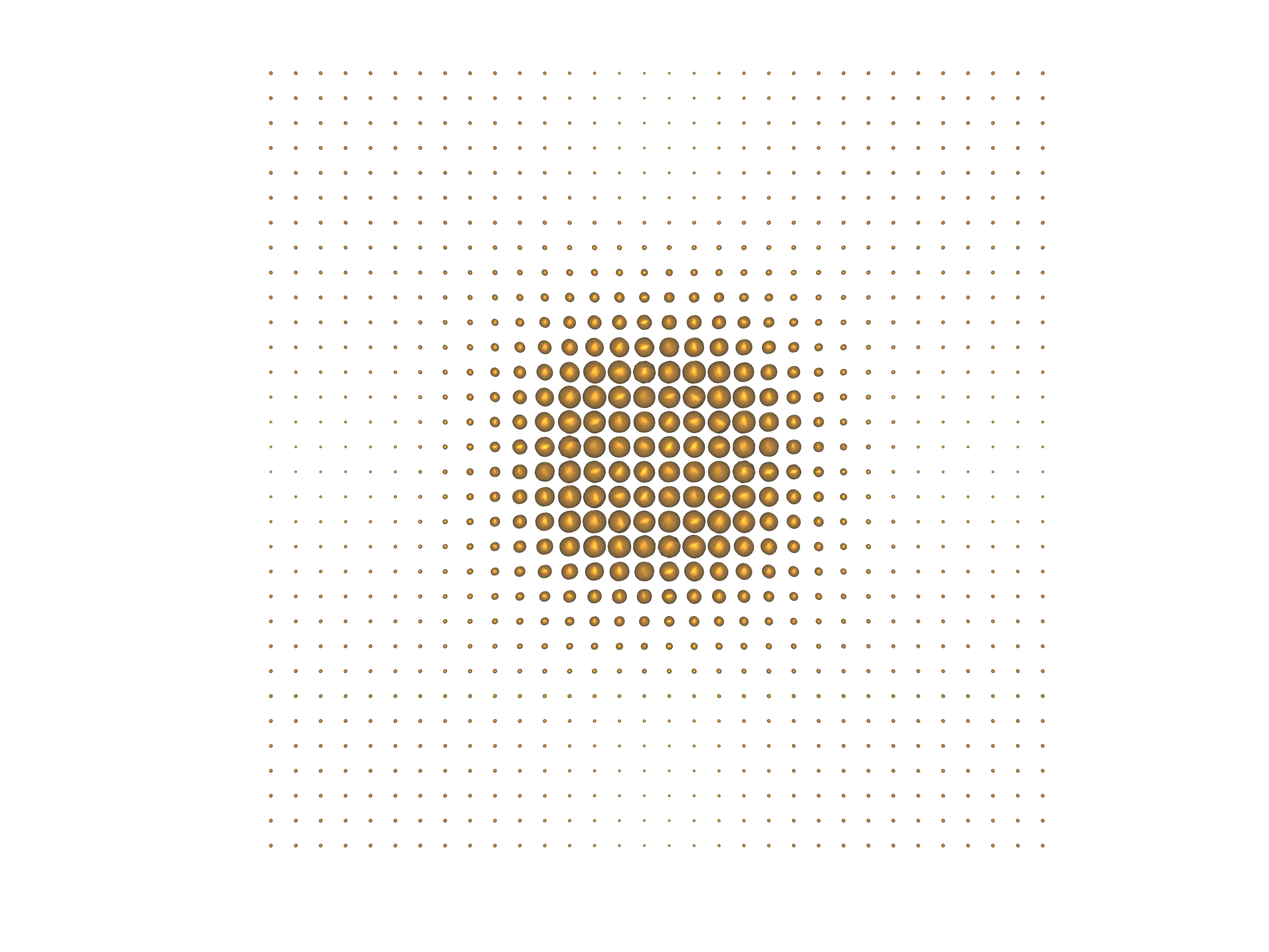}}
\subfloat[$t=0.2$]{\includegraphics[width=0.34\textwidth]{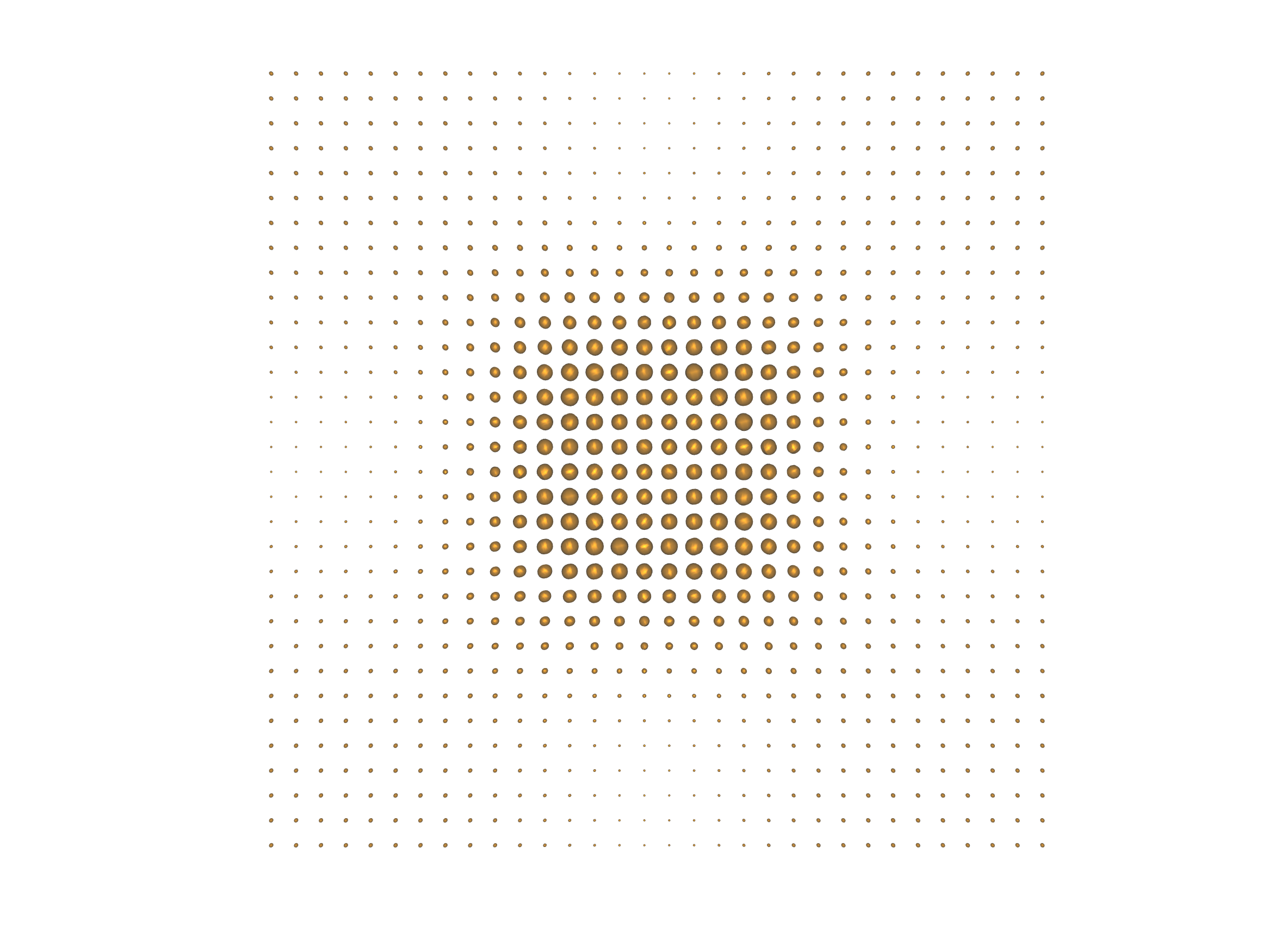}}
\subfloat[$t=0.3$]{\includegraphics[width=0.34\textwidth]{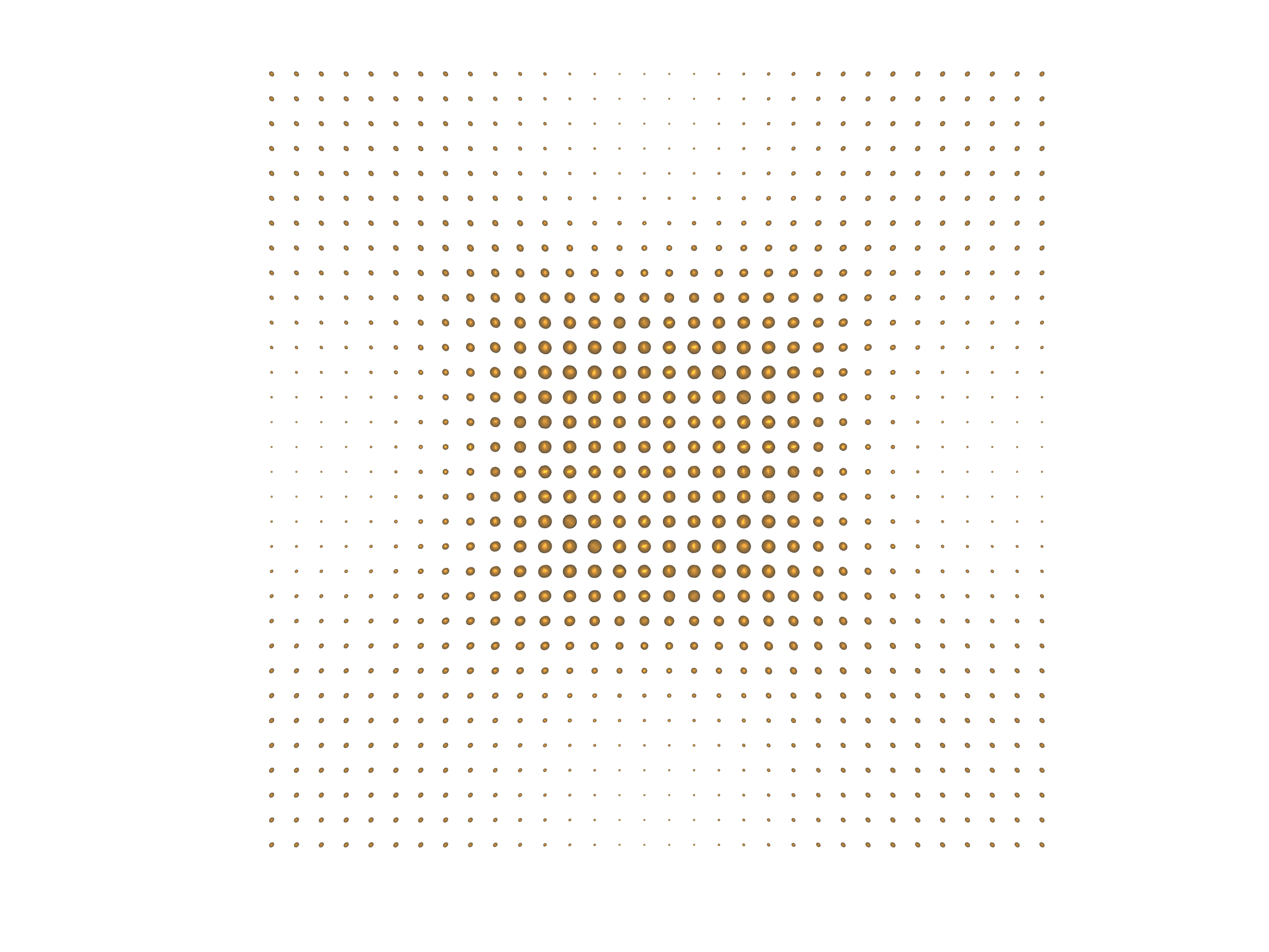}}
\\
\subfloat[$t=0.4$]{\includegraphics[width=0.34\textwidth]{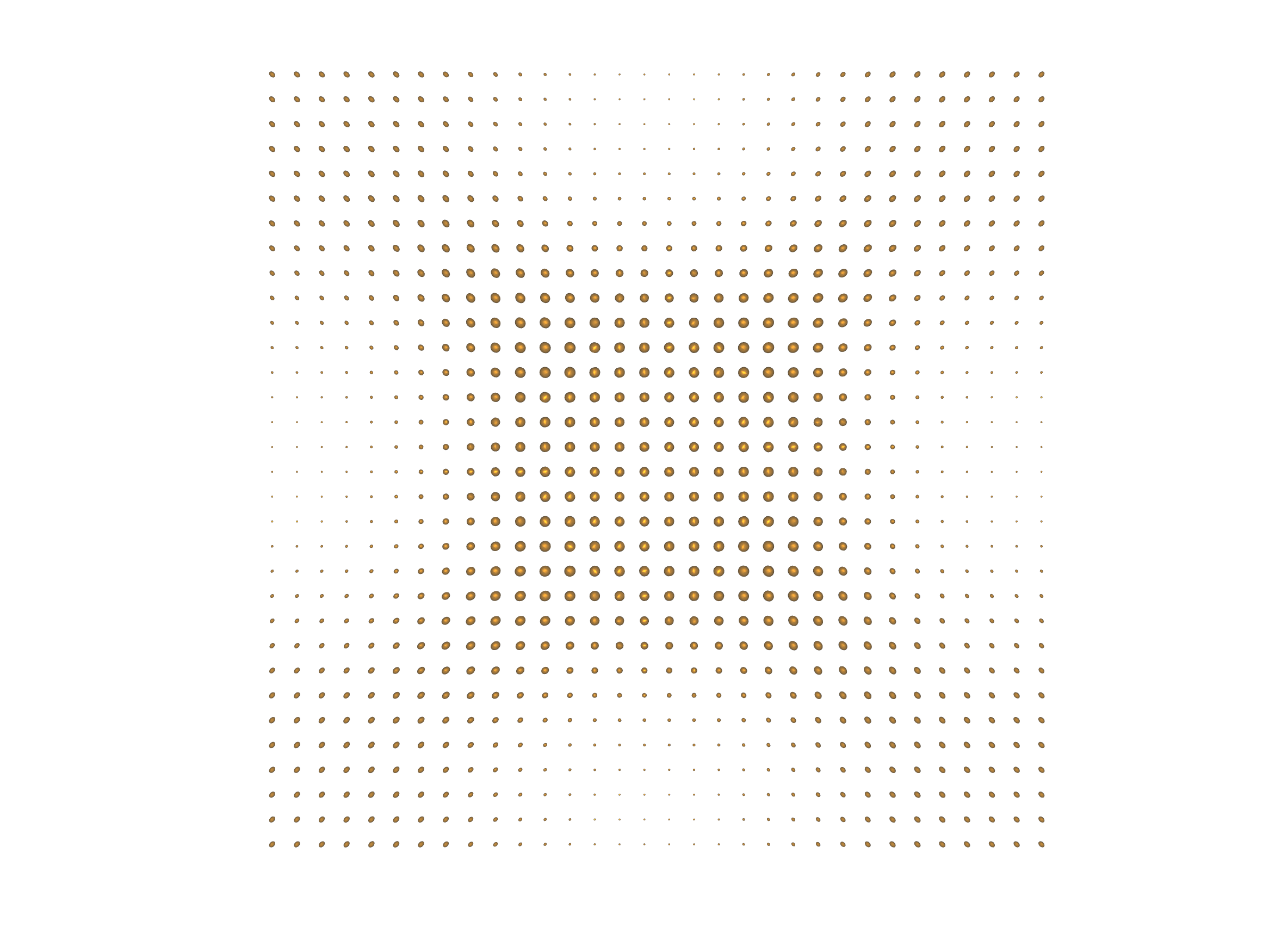}}
\subfloat[$t=0.5$]{\includegraphics[width=0.34\textwidth]{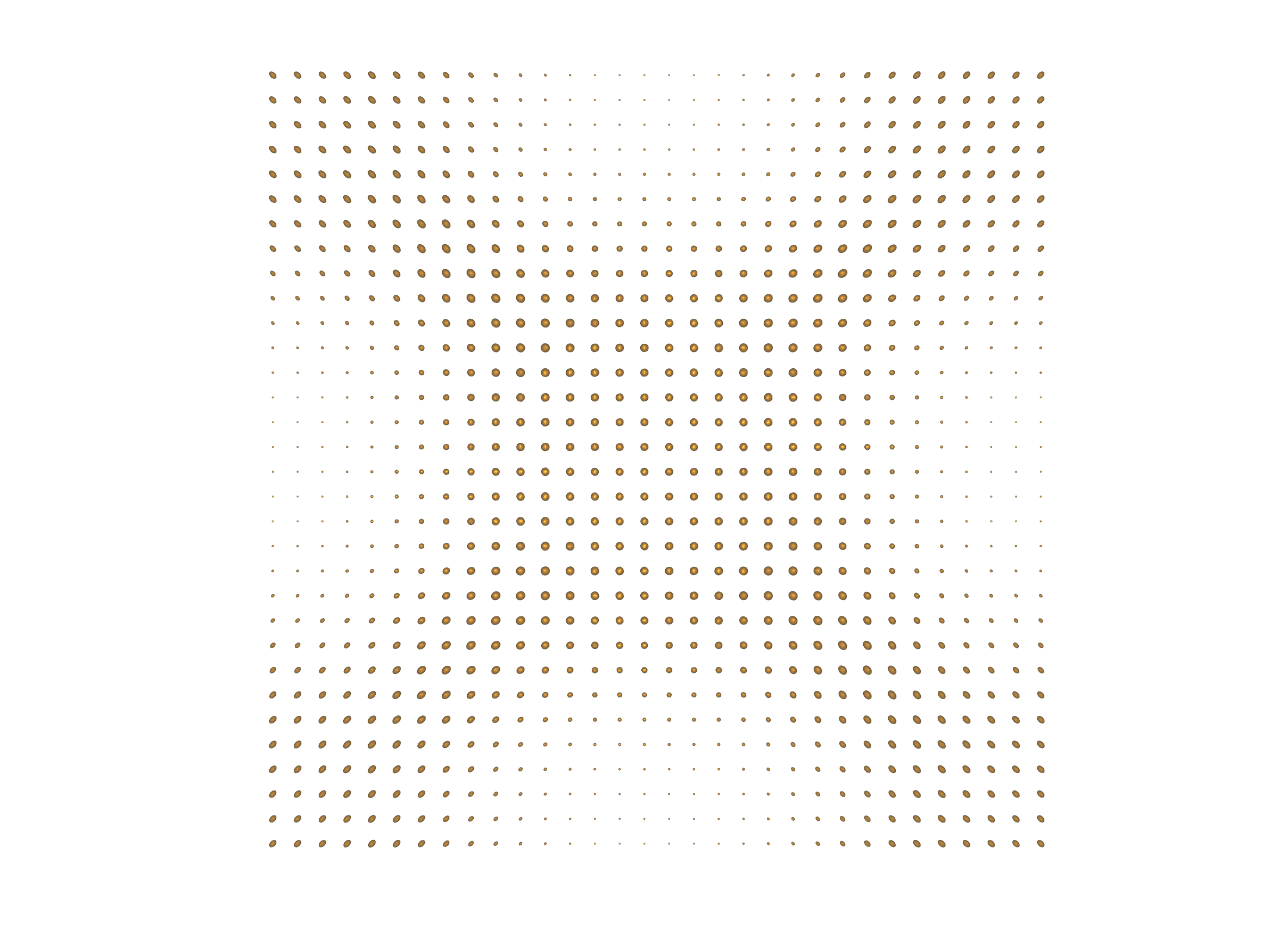}}
\subfloat[$t=0.6$]{\includegraphics[width=0.34\textwidth]{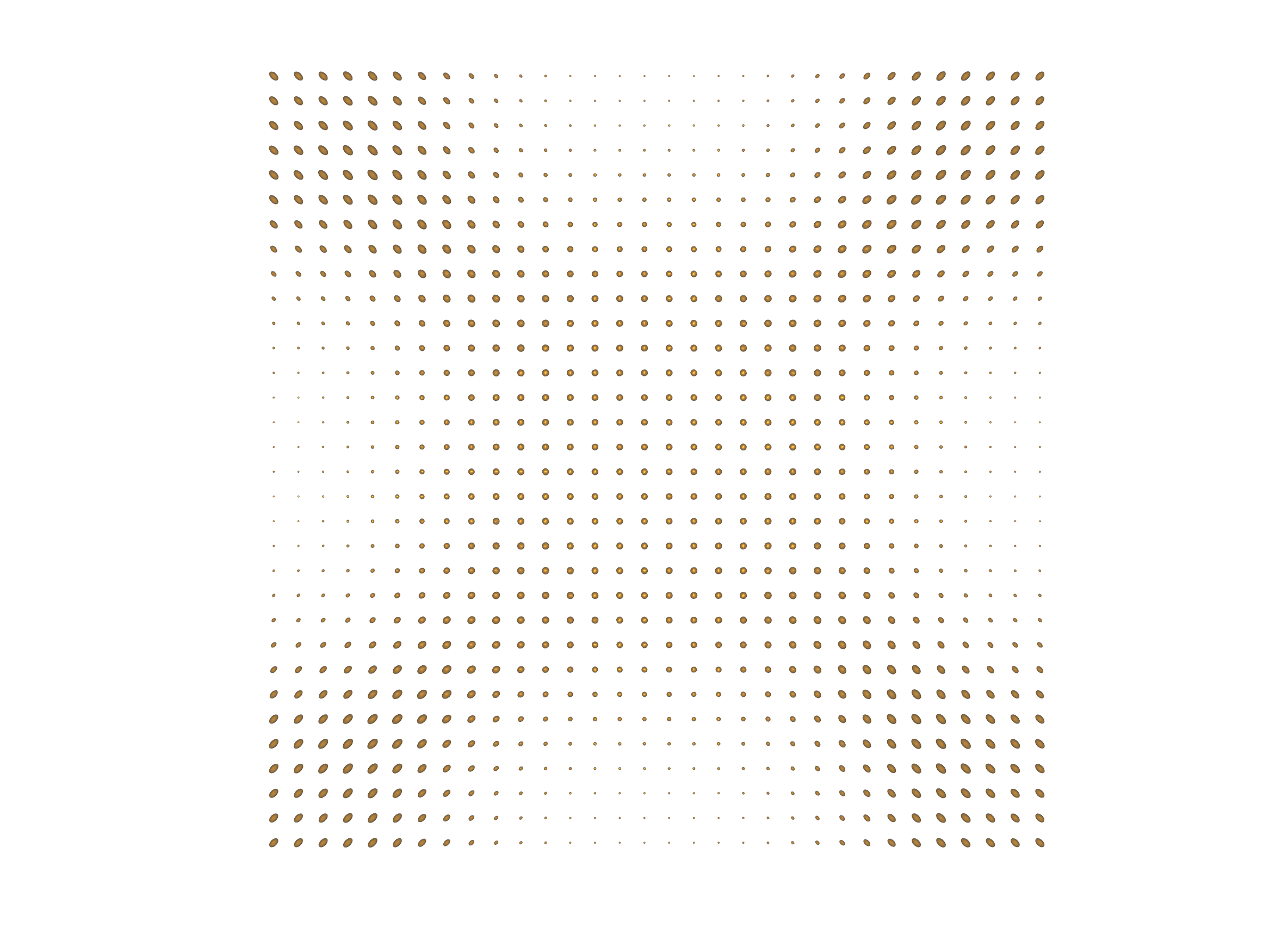}}
\\
\subfloat[$t=0.7$]{\includegraphics[width=0.34\textwidth]{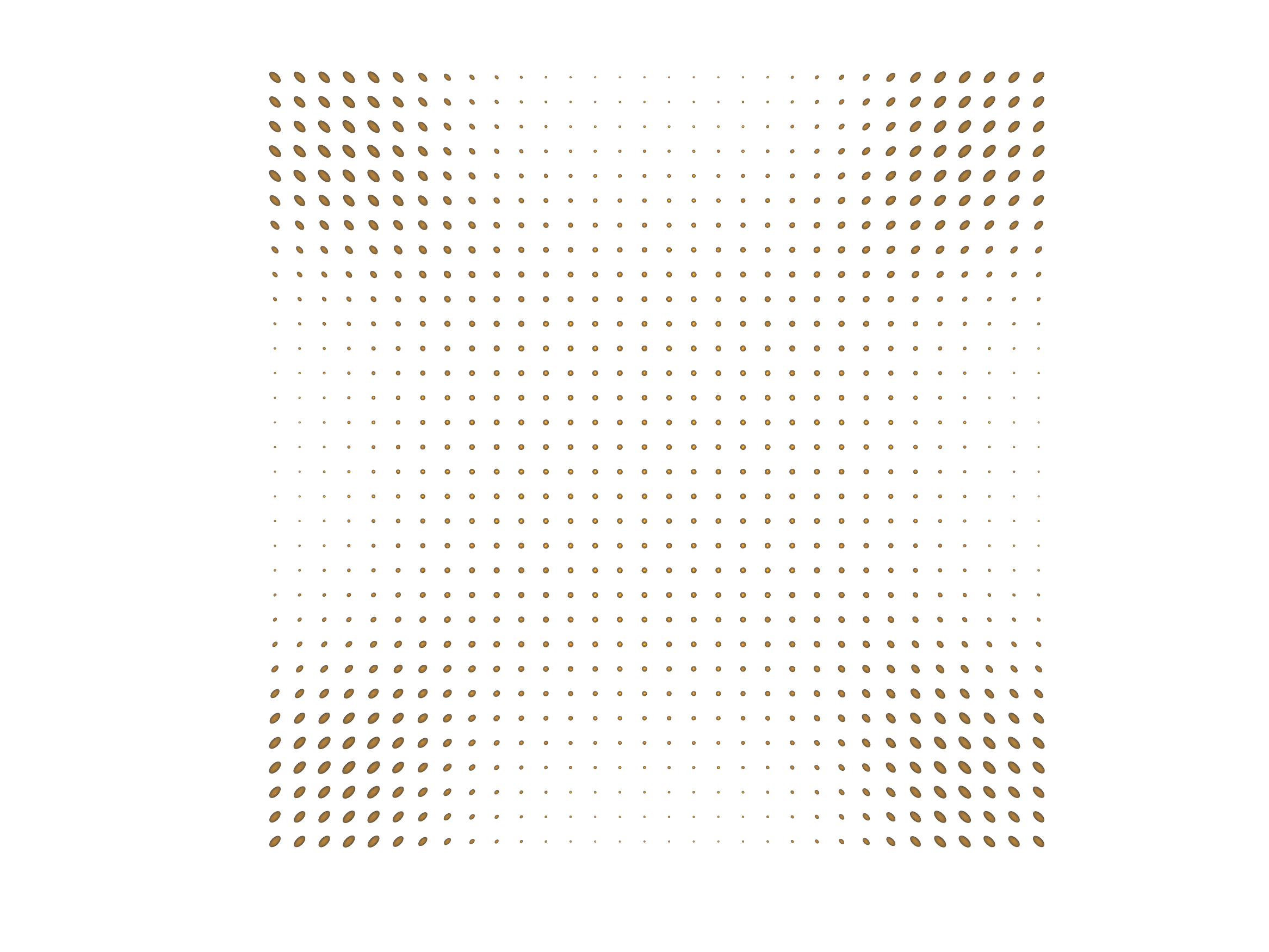}}
\subfloat[$t=0.8$]{\includegraphics[width=0.34\textwidth]{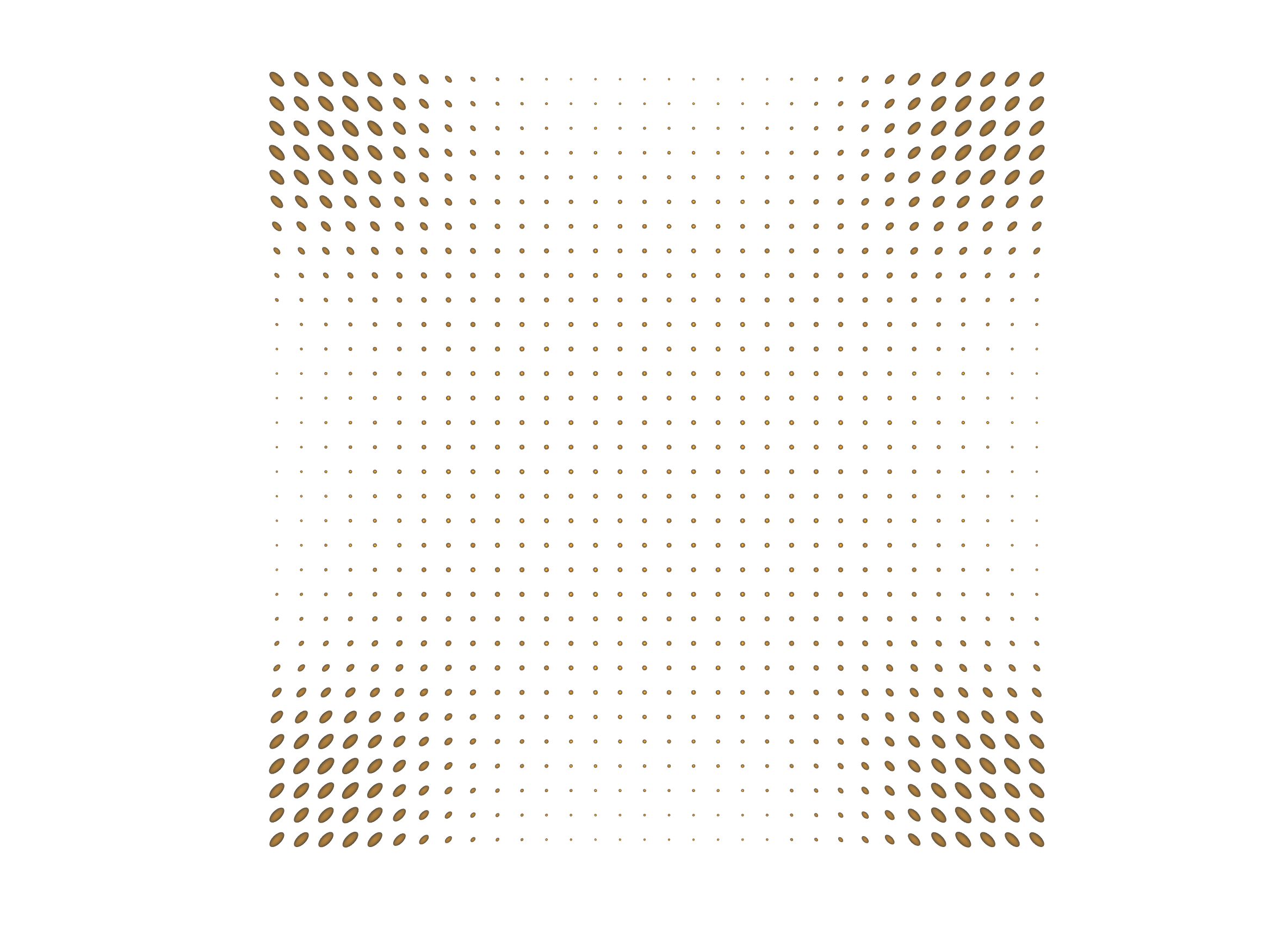}}
\subfloat[$t=0.9$]{\includegraphics[width=0.34\textwidth]{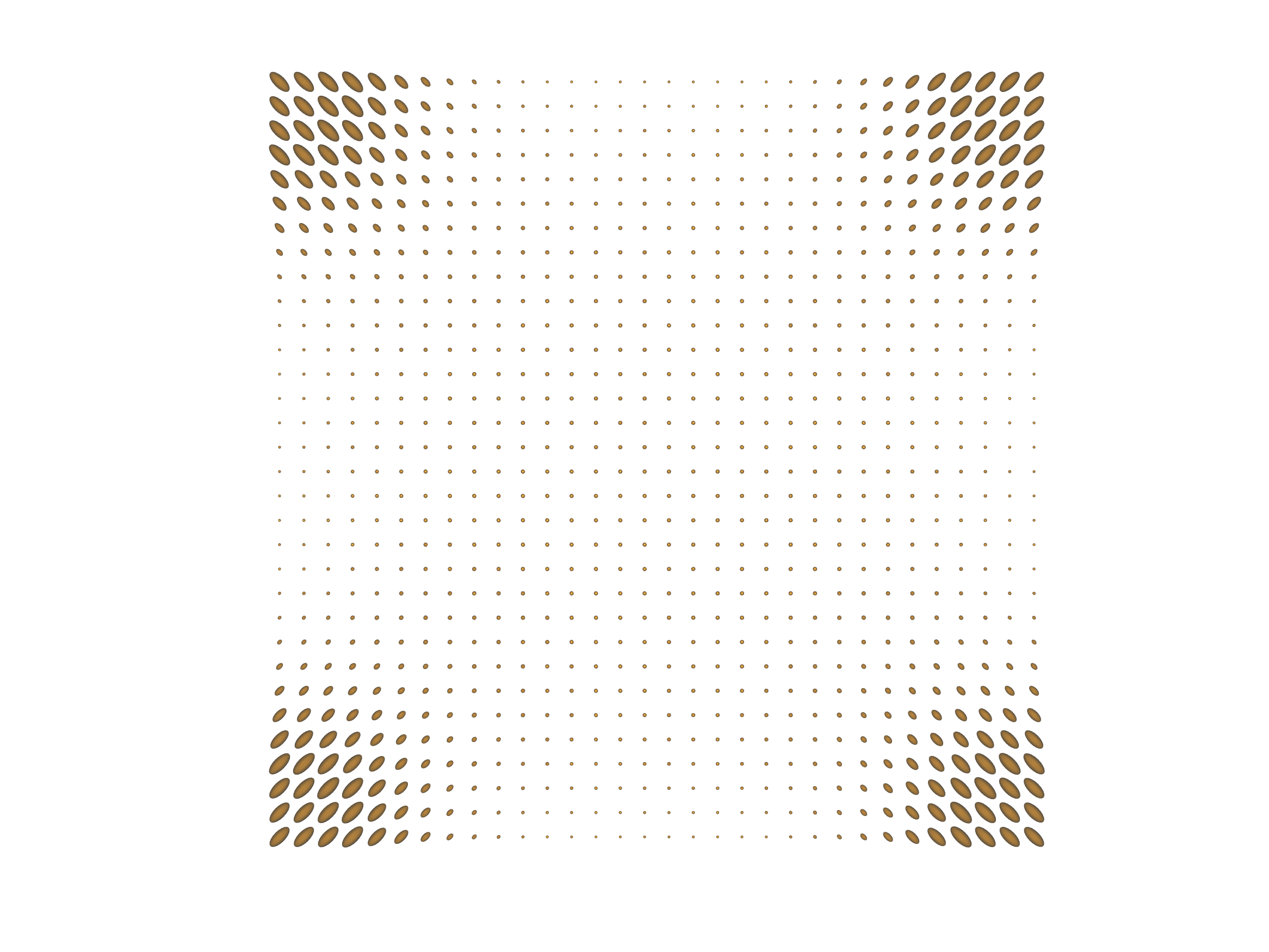}}
\caption{Interpolation with $\gamma=0.01$}
\label{fig:ex1interp}
\end{figure}

We tested our algorithm on a synthetic data set with $n=3$. The initial density is a disk positioned at the
center of the square domain and all the ellipsoids are isotropic.
The terminal density contains four quarter discs located at the corners of the square domain, and the
four components have different dominant directions. Both of them are depicted in
Figure~\ref{fig:ex1marginals}. The densities have been smoothed to have low density contrast $10$. Here the density contrast is defined to be the maximum of the ratios between the eigenvalues at different locations. In Figure~\ref{fig:ex1interp}, we show the optimal density flow with grid size $32\times 32 \times 10$ in space-time and parameter $\gamma=0.01$. The masses split into four components and the ellipsoids change gradually from isotropic to anisotropic.

To demonstrate the performance of our algorithm, we tested it on the same problem with different mesh grid sizes: $16 \times 16 \times 10$, $32\times 32 \times 20$, $64 \times 64 \times 40$ in space-time. We set the tolerance of the outer SQP iterations to $10^{-3},$ and that of the preconditioning conjugate gradient solver in each iteration to $10^{-3}$. The numbers of SQP iterations for convergence are shown in Table \ref{tab:matrixgrid1} for different mesh sizes.
\begin{table}
\centering
\begin{tabular}[h]{|| c | c ||}
\hline
Grid Size & SQP iterations \\
\hline
$16 \times 16\times 10$ & $19$ \\
$32 \times 32\times 20$ & $27$\\
$64 \times 64 \times 40$ & $35$\\
\hline
\end{tabular}
\caption{Number of SQP iterations required on different grid sizes for density contrast $10$.}
\label{tab:matrixgrid1}
\end{table}
\begin{table}
\centering
\begin{tabular}[h]{|| c | c ||}
\hline
Grid Size & SQP iterations \\
\hline
$16 \times 16\times 10$ & $25$ \\
$32 \times 32\times 20$ & $31$\\
$64 \times 64 \times 40$ & $62$\\
\hline
\end{tabular}
\caption{Number of SQP iterations required on different grid sizes for density contrast $50$.}
\label{tab:matrixgrid2}
\end{table}
\begin{table}
\centering
\begin{tabular}[h]{|| c | c ||}
\hline
Parameter $\gamma$ & SQP iterations \\
\hline
$1$ & $77$ \\
$0.1$ & $52$\\
$0.01$ & $31$\\
\hline
\end{tabular}
\caption{Number of SQP iterations required for different $\gamma$.}
\label{tab:matrixpara}
\end{table}

We then studed the influence of density contrast and the parameter $\gamma$ on the number of iterations needed to converge.
The results for density contrast $50$ are shown in Table~\ref{tab:matrixgrid2} with tolerance $10^{-2}$.
We can see that the number of iterations increases as we increase the density contrast.
Table~\ref{tab:matrixpara} showcases the results for different $\gamma$ values with fixed grid size $32 \times 32 \times 20$.
We observe that the number of iterations is positively correlated with the value of $\gamma$.

Finally, we test our algorithm on a 3D data set.
Table~\ref{tab:matrixgrid3} displays the number of iterations for different grid sizes with density contrast $30$ and parameter $\gamma=0.1$.
\begin{table}
\centering
\begin{tabular}[h]{|| c | c ||}
\hline
Grid Size & SQP iterations \\
\hline
$16 \times 16\times 16 \times 10$ & $19$ \\
$32 \times 32\times 32 \times 10$ & $25$\\
$64 \times 64 \times  64 \times 10$ & $23$\\
\hline
\end{tabular}
\caption{Number of SQP iterations required on different grid sizes for 3D densities.}
\label{tab:matrixgrid3}
\end{table}

\subsection{Vector case}
An important application of vector-valued optimal mass transport is color image processing.
In this cases, the vector-valued densities have three components corresponding to the intensities of the three basic colors red (R), green (G) and blue (B).
The masses can transfer from one color channel to another and the cost of transferring is captured using a weighted graph $\cF$.
Here, we treat the three colors equally and take the graph to be a complete graph with unit weights, namely, $W=I$ and
	\[
		\mD =
		\left[\begin{matrix}
		1 & 1 & 0\\
		-1 & 0 & 1\\
		0 & -1 & -1
		\end{matrix}\right].
	\]
The matrices $\mD_1,\mD_2$ in \eqref{eq:vecomtcvx} are then
	\[
		\mD_1 =
		\left[\begin{matrix}
		1 & 1 & 0\\
		0 & 0 & 1\\
		0 & 0 & 0
		\end{matrix}\right],\quad
		\mD_2 =
		\left[\begin{matrix}
		0 & 0 & 0\\
		1 & 0 & 0\\
		0 & 1 & 1
		\end{matrix}\right].
	\]
	
\begin{figure}[h]
\centering
\subfloat[$\rho^0$]{\includegraphics[width=0.40\textwidth]{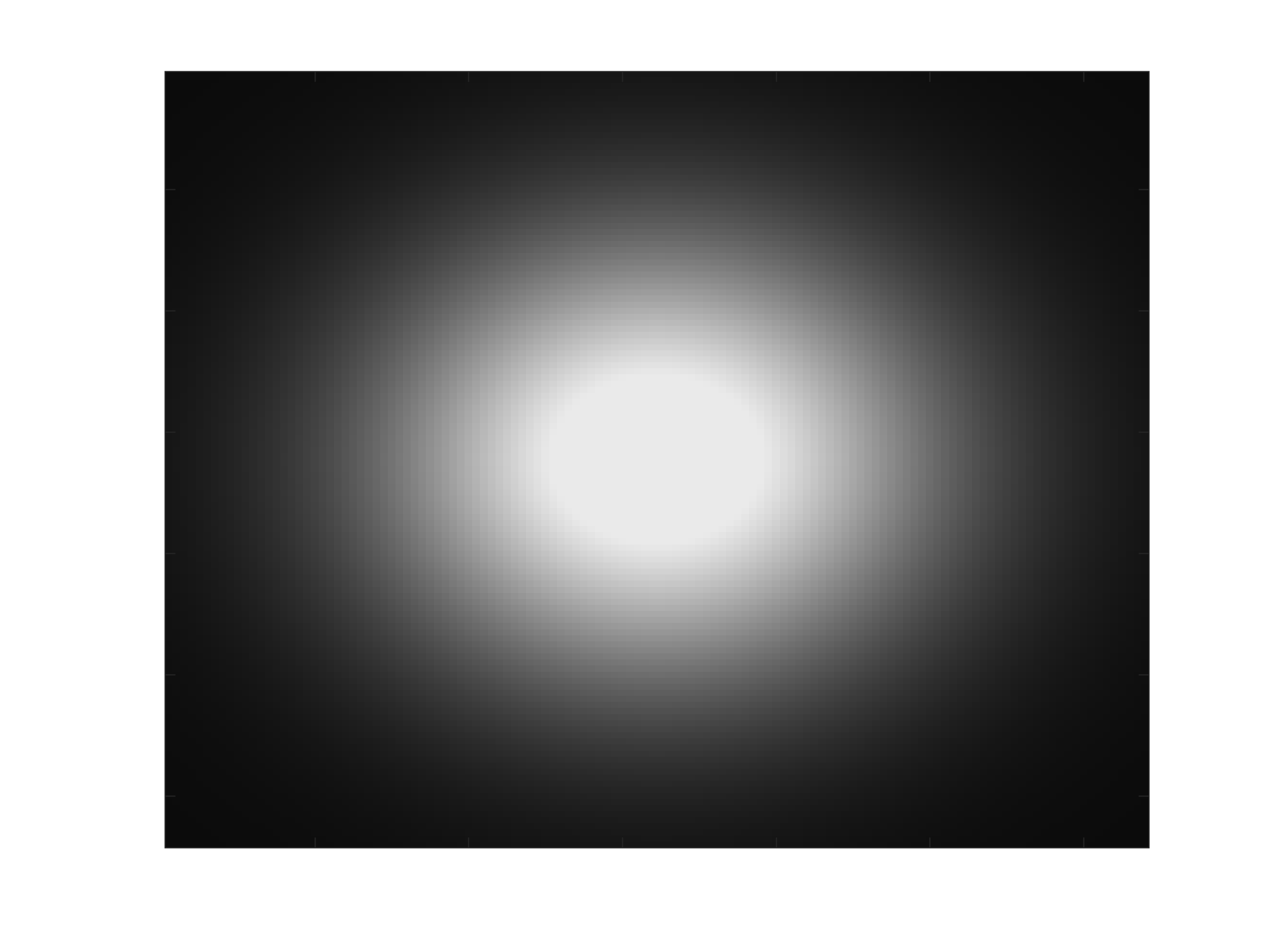}}
\subfloat[$\rho^1$]{\includegraphics[width=0.40\textwidth]{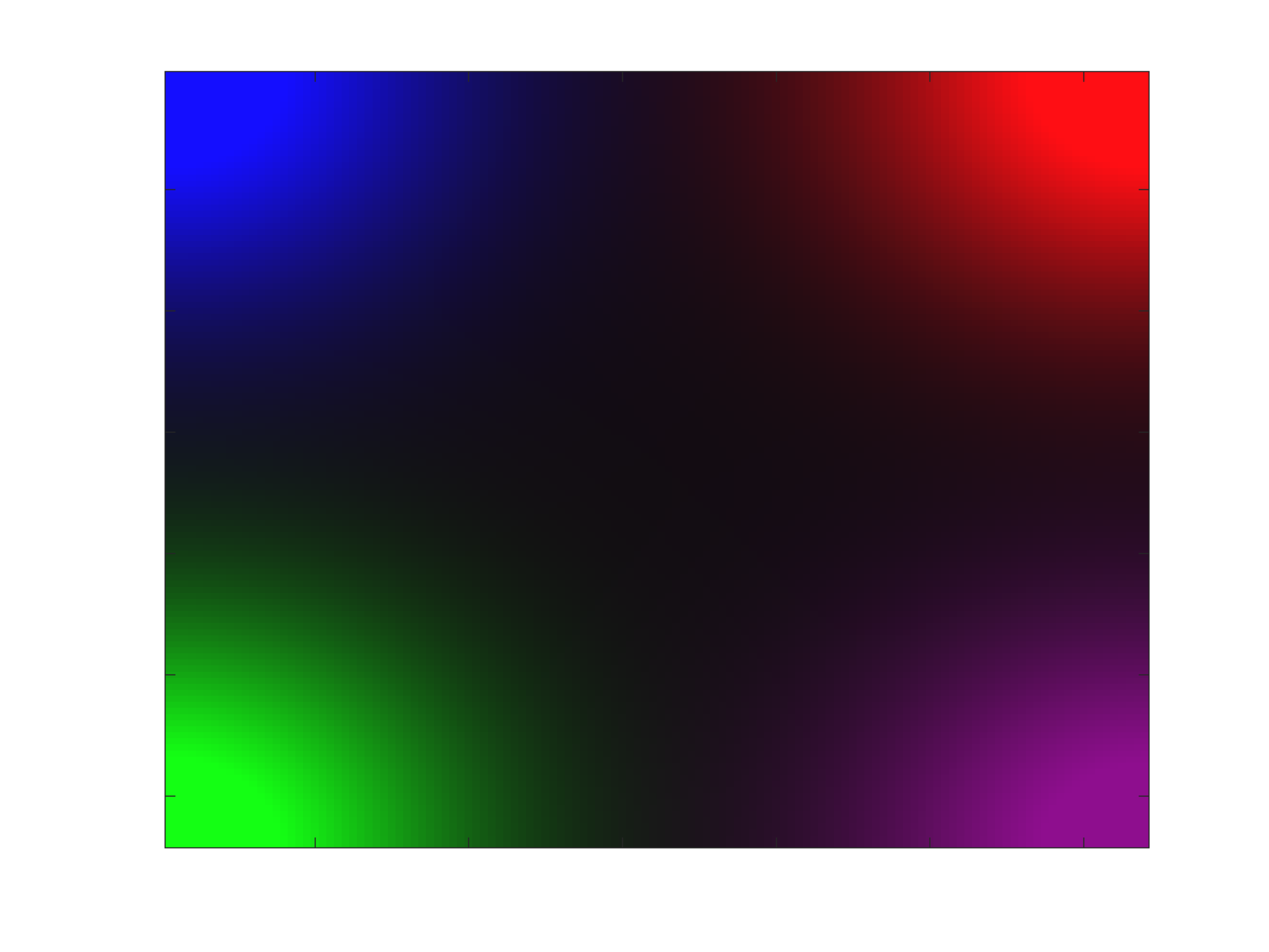}}
 \caption{Marginal distributions}
 \label{fig:ex2marginals}
\end{figure}

\begin{figure}[h]
\centering
\subfloat[$t=0.1$]{\includegraphics[width=0.2\textwidth]{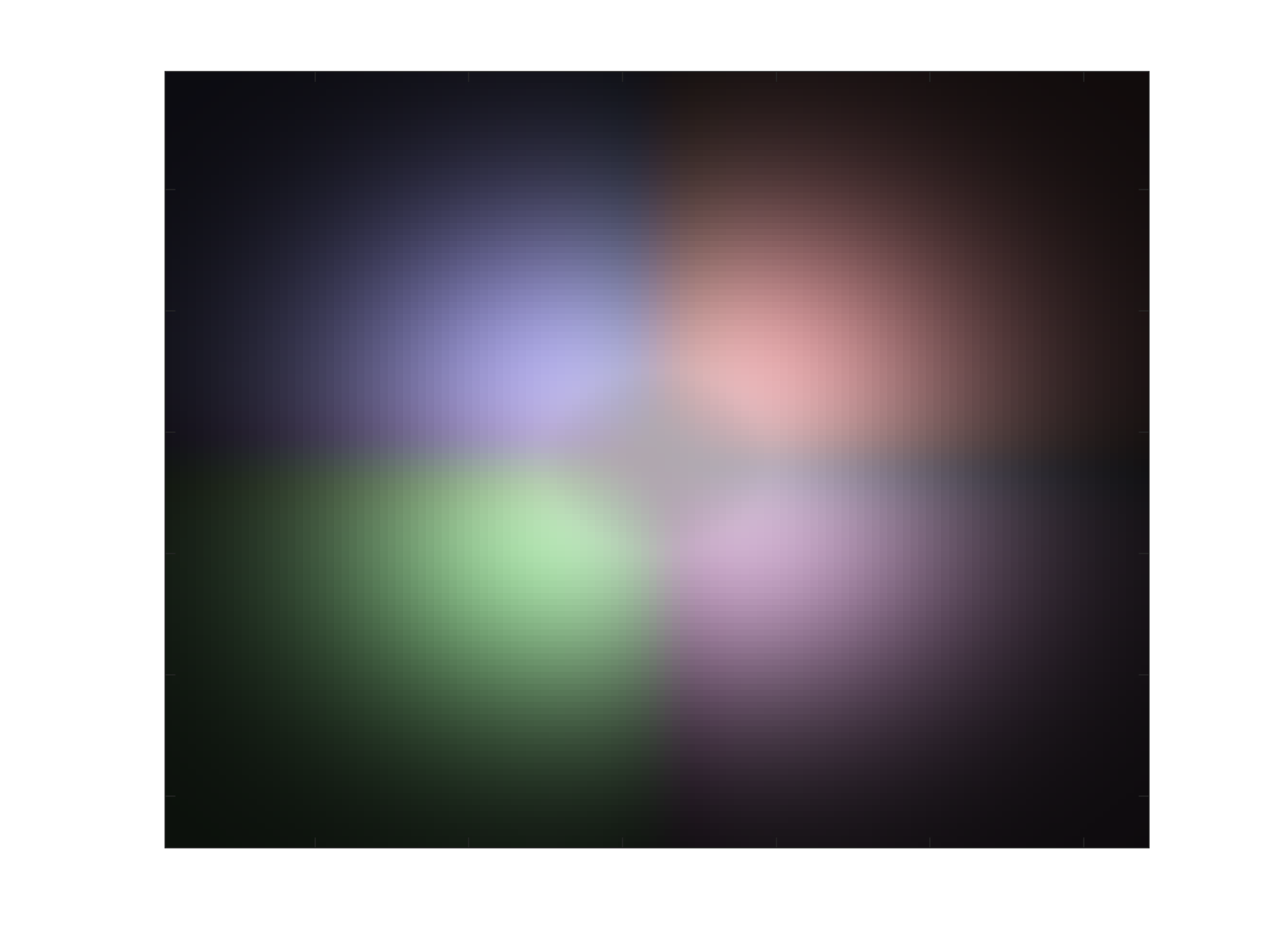}}
\subfloat[$t=0.2$]{\includegraphics[width=0.2\textwidth]{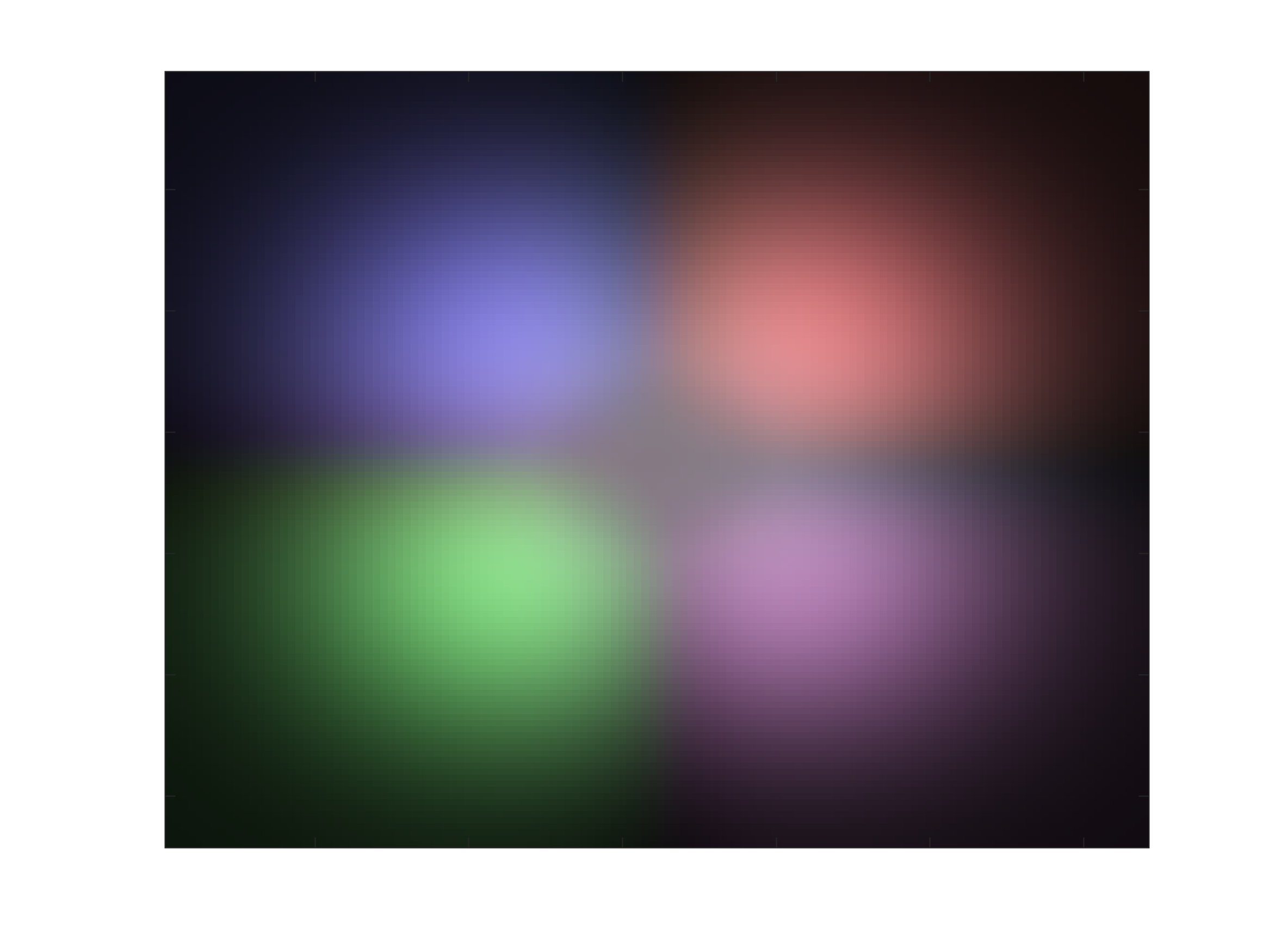}}
\subfloat[$t=0.3$]{\includegraphics[width=0.2\textwidth]{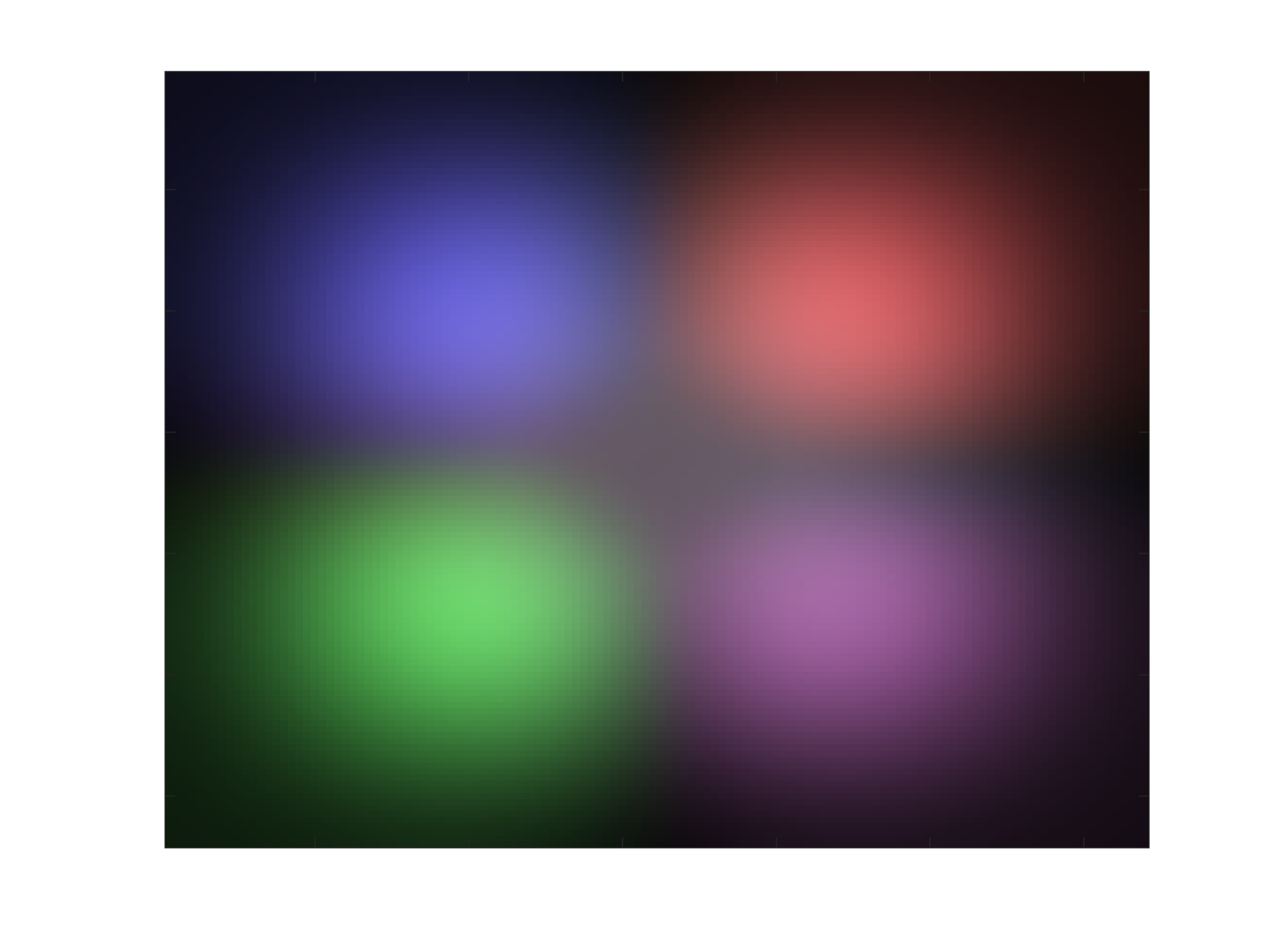}}
\subfloat[$t=0.4$]{\includegraphics[width=0.2\textwidth]{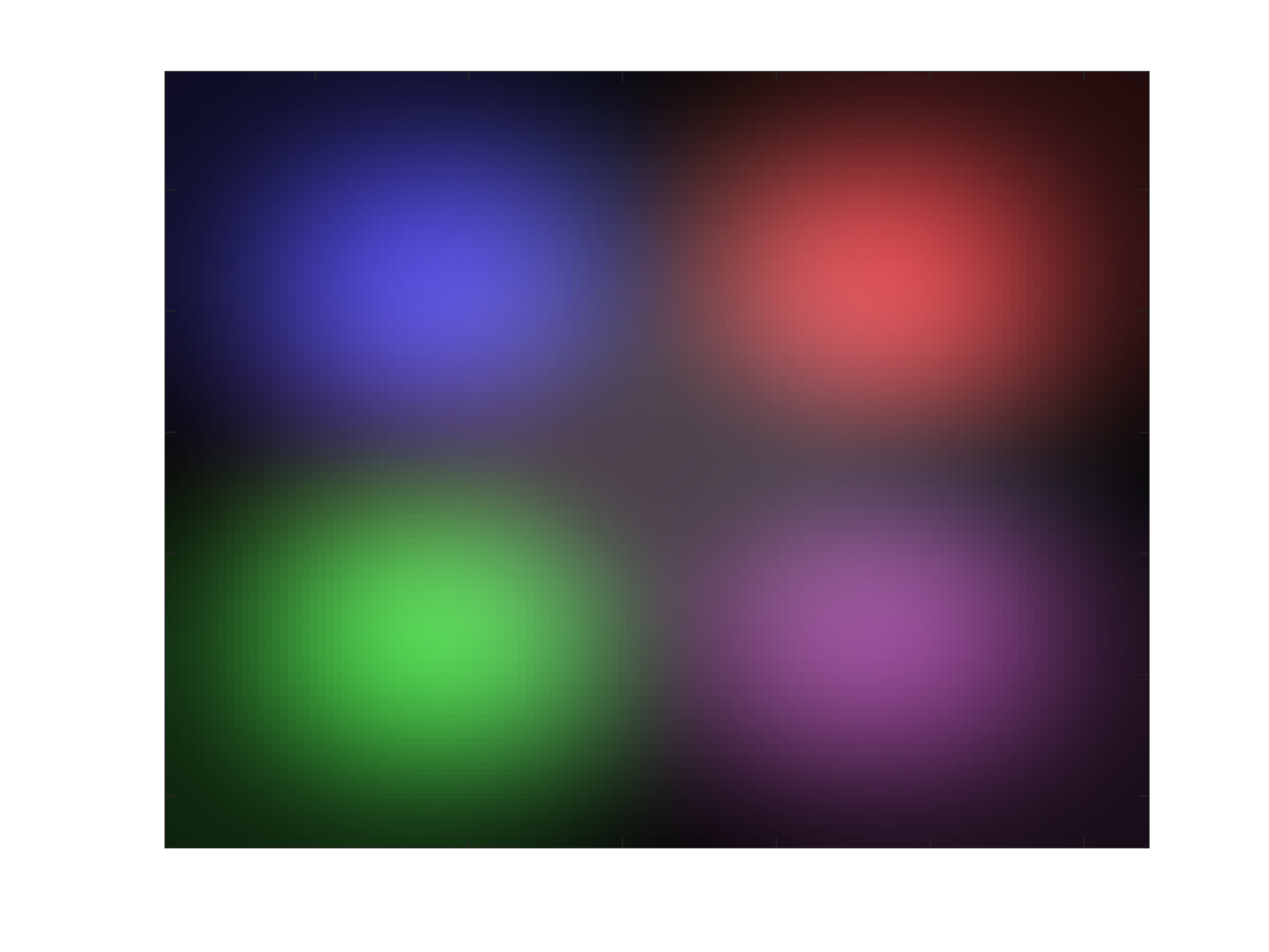}}
\\
\subfloat[$t=0.5$]{\includegraphics[width=0.2\textwidth]{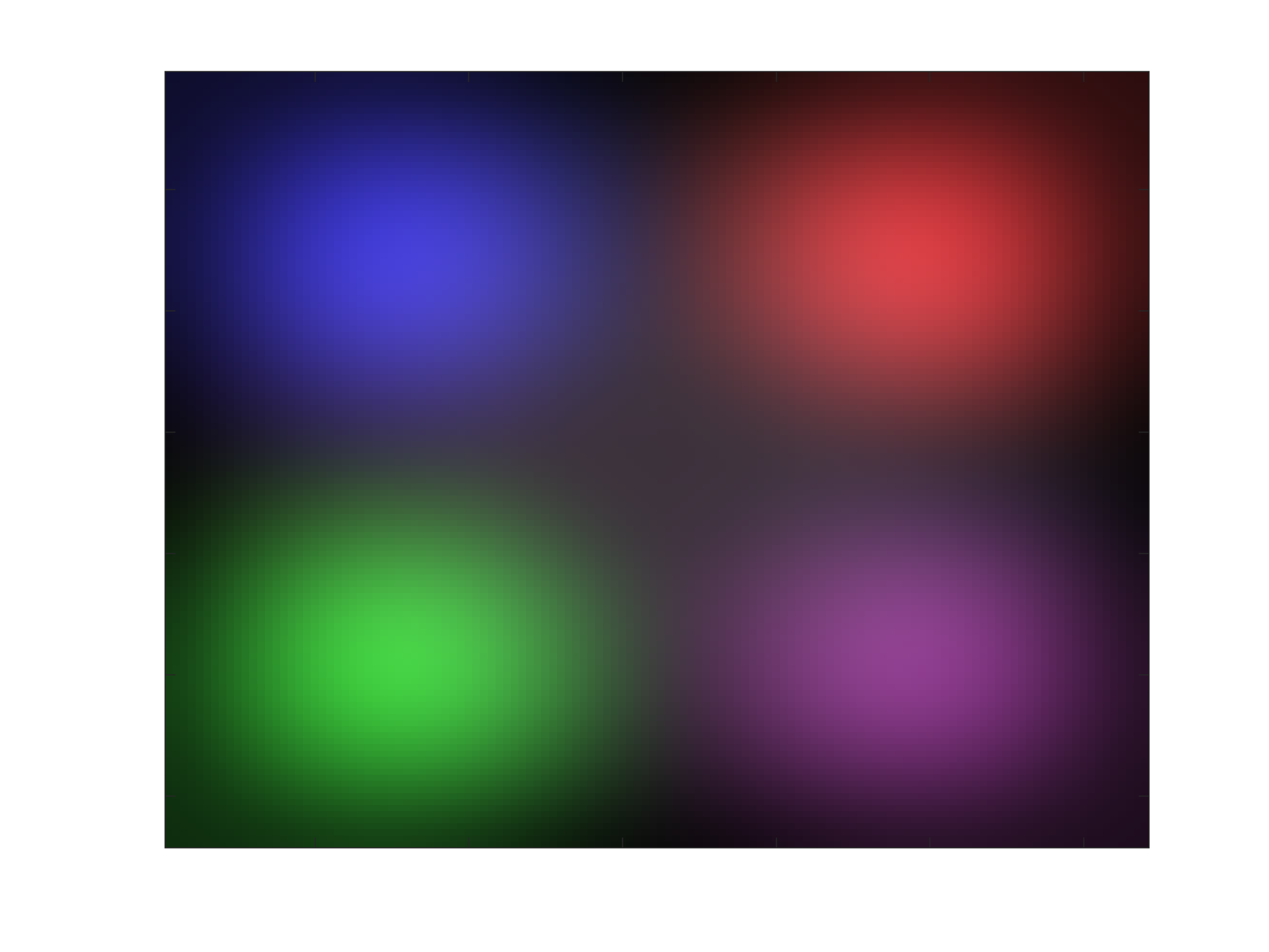}}
\subfloat[$t=0.6$]{\includegraphics[width=0.2\textwidth]{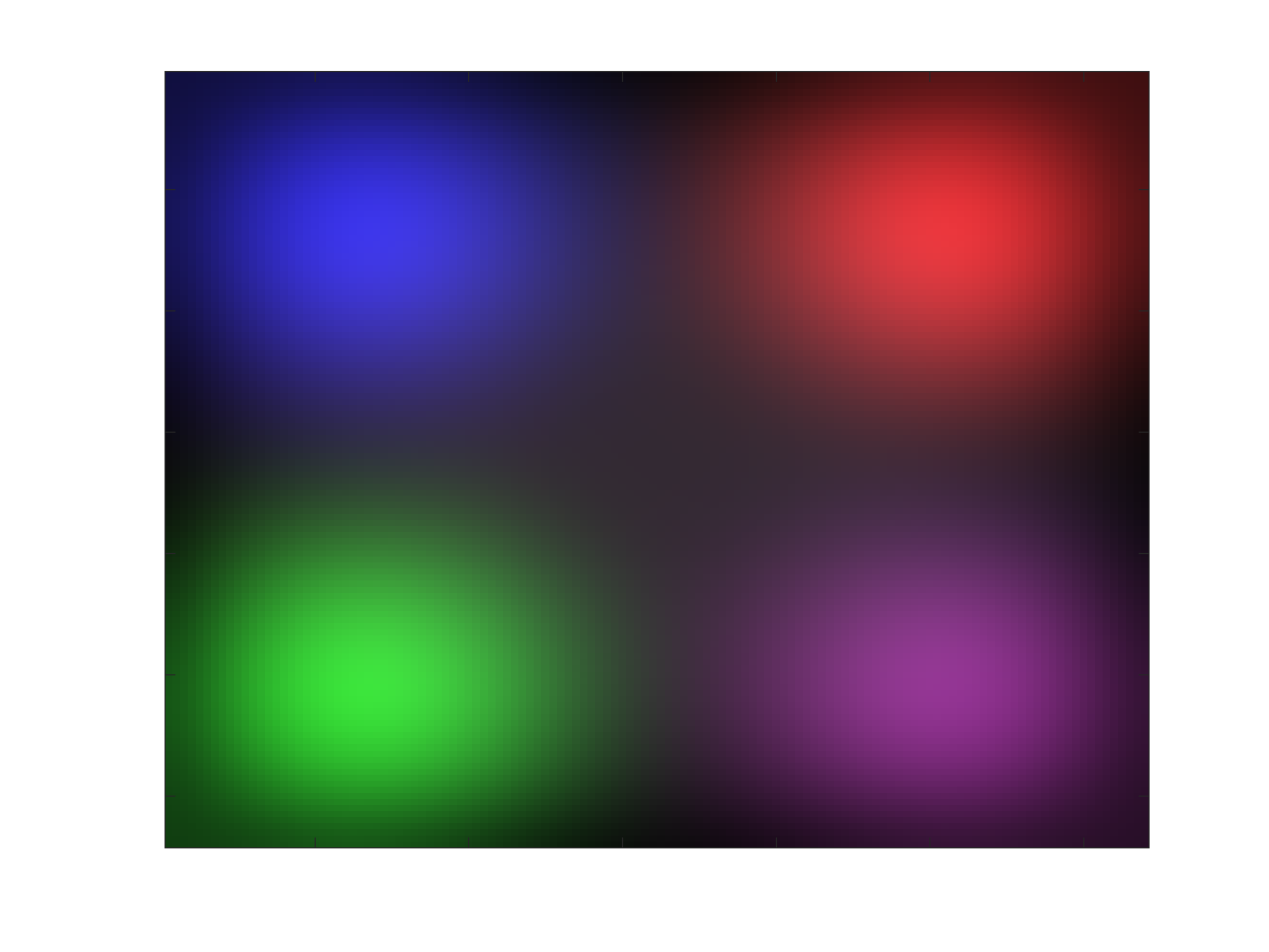}}
\subfloat[$t=0.7$]{\includegraphics[width=0.2\textwidth]{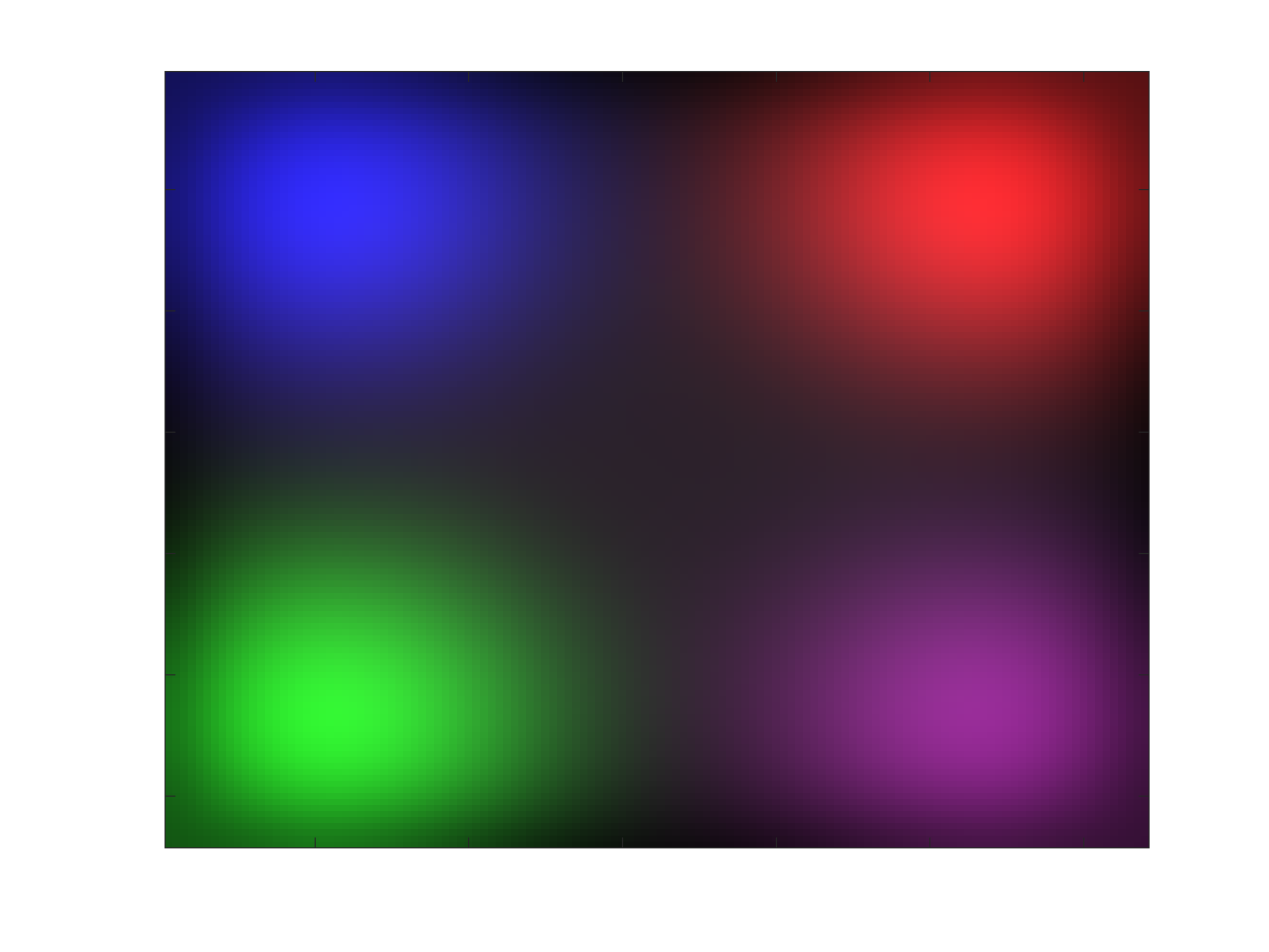}}
\subfloat[$t=0.8$]{\includegraphics[width=0.2\textwidth]{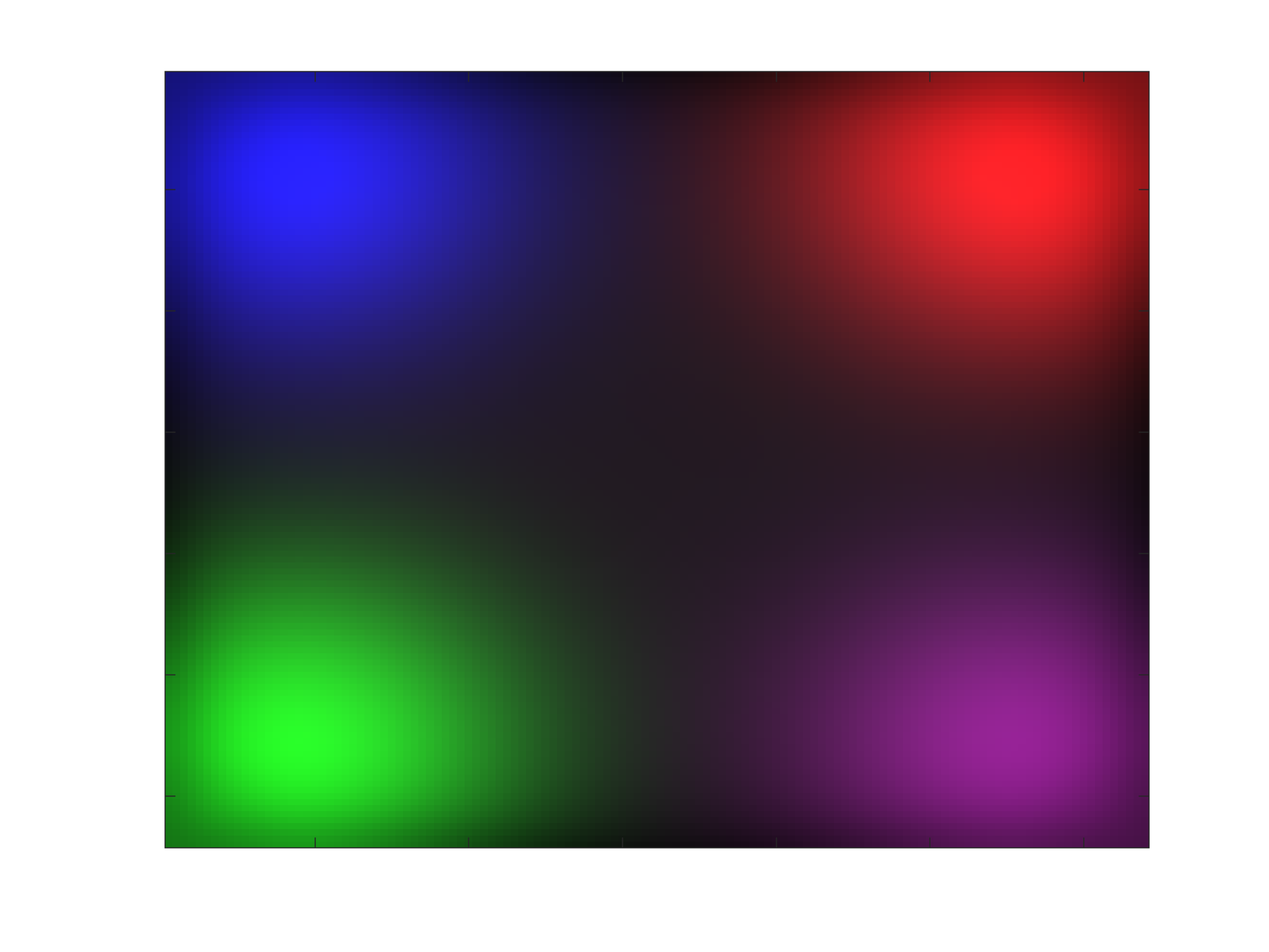}}
\subfloat[$t=0.9$]{\includegraphics[width=0.2\textwidth]{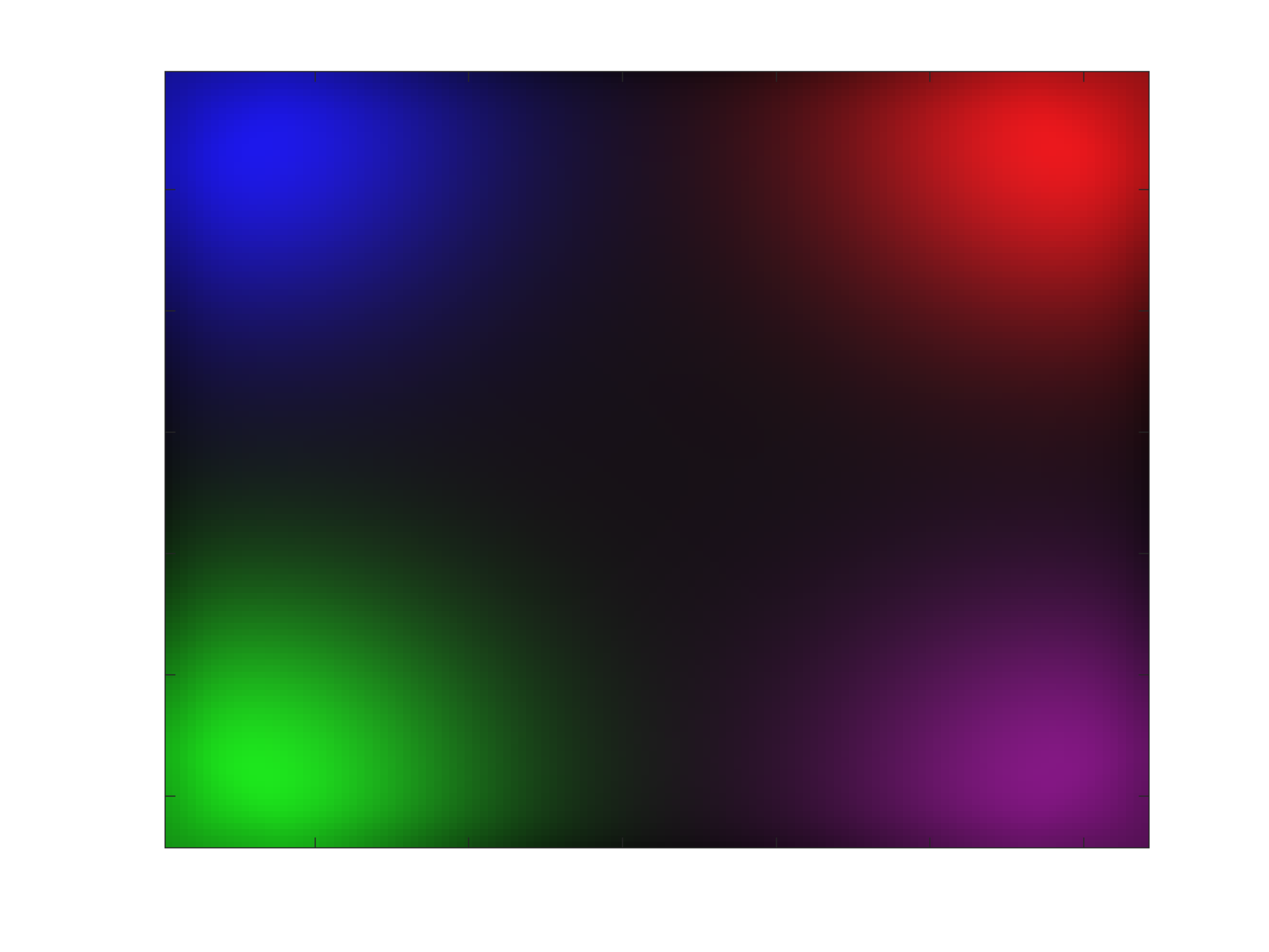}}
\caption{Interpolation with $\gamma=0.01$}
\label{fig:ex2interp}
\end{figure}

The two marginal densities are depicted in Figure~\ref{fig:ex2marginals}. The initial image $\rho^0$ is a disk located in the center of the square in white color, i.e.,
all three colors have equal intensity. The terminal distribution $\rho^1$ is an image of four circle quarters; one at each corner in different colors.
Both the images have been smoothed to have density contrast $\max_k \sup_{x,y} \rho_k^i(x)/\rho_k^i(y) \approx 10$. Figure~\ref{fig:ex2interp} illustrates
the optimal interpolation using vector-valued optimal transport with grid size $128\times 128\times 10$ in space-time and parameter $\gamma=0.01$.
We observe that the white disk split into four circle quarters and meanwhile the colors change gradually from white to four different colors.

We next tested the performance of the algorithm with respect to the grid size. For this,
we consider a grid hierarchy from a coarse grid of $32\times 32\times 10$ in space and time through a grid of $64 \times 64 \times 20$ to a grid of $128 \times 128 \times 40$.
The parameter $\gamma$ is set to be $0.01$. The tolerance for the outer SQP iteration is set to be $10^{-3}$ and in each iteration
the linear equation is solved with a relative residual of $10^{-2}$. The numbers of SQP iterations are recorded in Table~\ref{tab:vecgrid},
from which we observe that the number of iterations needed doesn't increase much as we increase the size of the mesh grids.
\begin{table}
\centering
\begin{tabular}[h]{|| c | c ||}
\hline
Grid Size & SQP iterations \\
\hline
$32 \times 32\times 10$ & $11$ \\
$64 \times 64\times 20$ & $12$\\
$128 \times 128 \times 40$ & $14$\\
\hline
\end{tabular}
\caption{Number of SQP iterations required on different grid sizes for density contrast $10$.}
\label{tab:vecgrid}
\end{table}

We also applied the same algorithm to images with a higher density contrast $100$. The results are shown in Table~\ref{tab:vecgrid1} for different grid sizes.
As can be seen from the table, increasing the density contrast leads to an increasing of the number of SQP iterations.
Again, the number of iterations needed to achieve certain precision is affected by the parameter. In Table~\ref{tab:vecpara} we display this
change as a function of $\gamma$ for fixed grid size $64\times 64\times 20$ and density contrast $100$.
\begin{table}
\centering
\begin{tabular}[h]{|| c | c ||}
\hline
Grid Size & SQP iterations \\
\hline
$32 \times 32\times 10$ & $24$ \\
$64 \times 64\times 20$ & $27$\\
$128 \times 128 \times 40$ & $32$\\
\hline
\end{tabular}
\caption{Number of SQP iterations required on different grid sizes for density contrast $100$.}
\label{tab:vecgrid1}
\end{table}
\begin{table}
\centering
\begin{tabular}[h]{|| c | c ||}
\hline
Parameter $\gamma$ & SQP iterations \\
\hline
$1$ & $48$ \\
$0.1$ & $42$\\
$0.01$ & $27$\\
\hline
\end{tabular}
\caption{Number of SQP iterations required for different $\gamma$.}
\label{tab:vecpara}
\end{table}

\section{Conclusions and future work} \label{sec:conclusions}

In this paper, we described a fast algorithm for the numerical implementation of both matrix-valued and vector-valued versions of optimal mass transport.
It is straightforward to extend this algorithm to cover matrix-valued transport problems with unequal masses (``unbalanced mass transport'') \cite{CheGeoTan17b}.
In the future, we intend to apply this methodology to various problems including diffusion tensor magnetic resonance data, biological networks,
and various types of vector-valued image data such as color and texture imagery. Finally, applying a multigrid methodology may speed up the
linear solver even further, and will be a future direction in our research.

\section*{Acknowledgements}
This project was supported by AFOSR grants (FA9550-15-1-0045 and FA9550-17-1-0435), grants from the National Center for Research Resources (P41-
RR-013218) and the National Institute of Biomedical Imaging and Bioengineering (P41-EB-015902), National Science Foundation (NSF), and
grants from National Institutes of Health (1U24CA18092401A1, R01-AG048769).

\bibliographystyle{IEEEtran}

\end{document}